\documentclass[11pt]{article}
\usepackage{geometry}
\usepackage[utf8]{inputenc}
\usepackage[english]{babel}
\usepackage[T1]{fontenc}
\usepackage{amsmath}
\usepackage{amsthm}
\usepackage{amsfonts}
\usepackage{amssymb}
\usepackage{bm}
\usepackage{lmodern}
\usepackage{verbatim}
\usepackage{ulem}

\usepackage{enumitem}

\usepackage[ruled]{algorithm}
\usepackage{algorithmic}

\usepackage{stmaryrd}

\usepackage{cancel}

\usepackage{hyperref}
\usepackage{url}

\usepackage{subfigure}
\usepackage{graphicx}

\usepackage{tikz}
\usepackage{tikz-3dplot}
\usetikzlibrary{3d}

\usetikzlibrary{positioning}
\usetikzlibrary{backgrounds}
\usetikzlibrary{patterns}
\usepackage{pgfplots}

\usepackage[font=small]{caption}

\usepackage{array}
\usepackage{graphbox}
\usepackage{xcolor,colortbl}

\usepackage{longtable}
\usepackage{booktabs}

\usepackage{wrapfig}
\usepackage{eurosym}
\usepackage{color}

\definecolor{Red}{rgb}{1,0,0}
\definecolor{Green}{rgb}{0,.6,0}
\definecolor{Blue}{rgb}{0,0,1}

\newcommand{\bl}{\color{black}} 

\usepackage{minted}

\usepackage{pifont}

\newcommand{\cmark}{\ding{51}}

\usepackage[title]{appendix}

\theoremstyle{definition}
\newtheorem{definition}{Definition}[section]

\theoremstyle{remark}
\newtheorem{remark}{Remark}[section]

\theoremstyle{plain}

\theoremstyle{plain}

\theoremstyle{plain}
\newtheorem{proposition}{Proposition}[section]

\theoremstyle{plain}
\newtheorem{assumption}{Assumption}[section]

\providecommand{\keywords}[1]
{
  \small
  \textbf{Key words.} #1
}

\providecommand{\jelcodes}[1]
{
  \small
  \textbf{MSC codes.} #1
}

\title{MOCVXPY: a CVXPY extension for multiobjective optimization}

\author{
  \href{mailto:ludovic.salomon@polymtl.ca}{\bl Ludovic Salomon}\thanks{
    GERAD and D\'epartement de math\'ematiques et de g\'enie industriel,
    \href{https://polymtl.ca}{\bl Polytechnique Montr\'eal},
    C.P. 6079, Succ. Centre-ville,
    Montr\'eal, Qu\'ebec, Canada H3C~3A7.
    \href{mailto:ludovic.salomon@polymtl.ca}{\bl \url{ludovic.salomon@polymtl.ca}},
  }
  \and
  \href{mailto:daniel.doerfler@uni-jena.de}{\bl Daniel Dörfler}\thanks{
    Department of applied mathematics,
    \href{https://www.uni-jena.de}{\bl Friedrich Schiller University Jena},
    Jena, Thuringia, Germany.
    \href{mailto:daniel.doerfler@uni-jena.de}{\bl \url{daniel.doerfler@uni-jena.de}},
  }
  \and
  \href{mailto:andreas.loehne@uni-jena.de}{\bl Andreas Löhne}\thanks{
    Department of applied mathematics,
    \href{https://www.uni-jena.de}{\bl Friedrich Schiller University Jena},
    Jena, Thuringia, Germany.
    \href{mailto:andreas.loehne@uni-jena.de}{\bl \url{andreas.loehne@uni-jena.de}},
  }
}

\date{August 2025}

\pgfplotsset{compat=1.18}

\begin{document}

\maketitle

\begin{abstract}
    MOCVXPY is an open-source Python library for convex vector optimization. It is built on top of CVXPY, a domain-specific language for single-objective convex optimization. MOCVXPY enables practitioners to describe their convex vector optimization problem in an intuitive algebraic language, that closely follows the mathematical formulation. This work presents the main features of MOCVXPY, explains some background of the algorithms it employs to solve the optimization problems, and illustrates its functionality through examples and two real-world applications in finance and energy. MOCVXPY is available at \url{https://github.com/salomonl/mocvxpy} under the Apache 2.0 licence, with some documentation and examples.
\end{abstract}


\jelcodes{90C29; 90C25; 90C59; 90B50; 97N80.}

\section{Introduction}\label{section:Introduction}

Multiobjective optimization (MO) is a generalization of classical single-objective optimization in which multiple objective functions must be optimized simultaneously~\cite{Collette2011MultiobjectiveOptimizationPrinciplesCaseStudies, Ehrgott2005MulticriteriaOptimization, Miettinen1999NonlinearMultiobjectiveOptimization}. In general, these objectives conflict with each other: improving one of them worsens some of the others. Consequently, the set of optimal solutions for such a problem is not a singleton (or a set of optimal decision vectors with the same objective values). It is composed of different feasible solutions, the so-called \textit{efficient set}. The values of these solutions in the objective space represent the best trade-offs between the different criteria and are referred to as the \textit{non-dominated set/Pareto front}. Once the Pareto front is obtained, decision makers can conscientiously select one of the generated optimal policies for a specific situation. Many applications fit into the framework: e.g., chemical engineering~\cite{Sharma2013}, computational bioinformatics~\cite{Handl2007}, machine learning~\cite{Alexandropoulos2019} or energy~\cite{Cui2017}.

Vector optimization (VO) refers to the optimization of a vector-valued objective function with respect to a partial order relation, determined by a suitable ordering cone \(\mathcal{C}\)~\cite{Ehrgott2005MulticriteriaOptimization, Jahn2009VectorOptimization, Lohne2011Vector}. In MO, the ordering cone is the non-negative orthant, making it a special case of VO. In the context of VO, the image of the set of efficient solutions by the objective function is called the set of minimal objective vectors with respect to \(\mathcal{C}\). Although less common than applications for MO, applications for VO can be found for example in financial mathematics~\cite{Ararat2023a, Feinstein2017, Hamel2013, Lohne2014a}, economics~\cite{Rudloff2021, Ruszczynski2003}, game theory~\cite{Feinstein2024a}, polyhedral calculus~\cite{Ciripoi2019} or medical imaging~\cite{Eichfelder2014}.

In general, it is impossible to find the entire non-dominated set/set of minimal objective vectors since it can be infinite and does not need to have a finite representation.
Existing algorithms aim to generate a discrete representation of this set that is ``sufficiently good'' to approximate the true solution set in the objective space.

Heuristics such as evolutionary~\cite{Branke2008MultiobjectiveOptimizationInteractiveEvolutionaryApproaches} or particle-swarm~\cite{Shami2022} methods are commonly used in the context of MO. These procedures evolve a ``population'' of points in the decision space via mutation/crossover operations/position updates that bring them closer to the non-dominated set along the iterations. Such methods are not exact in the sense that they do not have convergence guarantees. They are stochastic by nature and require a considerable number of function evaluations to deploy, which may make them prohibitive for large-scale problems. In general, they cannot exploit the analytical structure of the problem (e.g., convexity, differentiability) when it is available (some exceptions exist, e.g.,~\cite{Schutze2016a}, that combines a continuation method with a MO evolutionary algorithm). Their efficiency depends on the choice of their numerous parameters. On the other hand, the minimal assumptions they require make them suitable, if not effective, for a large number of applications~\cite{Sharma2022}. Practitioners have also access to many tested libraries that offer a vast catalog of such methods: e.g., \texttt{jMetal}~\cite{Durillo2011}, \texttt{pagmo}~\cite{Biscani2020}, \texttt{pymoo}~\cite{Blank2020}, \texttt{Platypus}~\cite{Hadka2024}, \texttt{DEAP}~\cite{Fortin2012}, \texttt{inspyred}~\cite{Garett2012}; which may favor their choice over potentially better-tailored procedures for a certain class of problems.

What should a practitioner use if its VO or MO problem possesses some analytical properties (e.g., derivatives, convexity)?
Scalarization-based algorithms reformulate the original MO/VO problem into a sequence of parameterized single-objective subproblems (see~\cite{Wiecek2016} for a comprehensive survey on this topic). By varying the parameter values, one can obtain different solutions of the initial problem. Since the Normal Boundary Intersection method was described in 1998~\cite{Das1998}, smart parameter choices that exploit the structure of the original problem have been proposed (e.g.,~\cite{Burachik2017, Burachik2022, Eichfelder2009}). However, these methods only offer cardinality guarantees. One can only assess that the discrete representation generated by these methods may have a maximum number of solutions depending on the number of combinations of parameter values involved in the optimization process. Note that other techniques can be used such as descent-based algorithms for MO and VO: see~\cite{Eichfelder2021b, Fukuda2014} for an overview.

Enclosure-based approaches (e.g.,~\cite{DeSantis2021, DeSantis2020, DeSantis2024, Eichfelder2021a,  Eichfelder2024a, Eichfelder2021, Eichfelder2023a, Eichfelder2024, Niebling2019}) and outer approximation approaches (also known as Benson-type or sandwiching algorithms), e.g., ~\cite{Ararat2023, Ararat2022, Ararat2024, Benson1998a, Dorfler2022, Ehrgott2012, Ehrgott2011, Helfrich2024a, Keskin2023,  Lammel2025, Lohne2014}, are a class of scalarization-based approaches that eliminate this limitation. They enclose the non-dominated set of solutions between an outer and inner approximation. Enclosure-based approaches use local upper/lower bound sets while outer approximation algorithms require convex polyhedra. This ``enclosure'' becomes finer along the iterations. When the ``distance'' between the outer and inner approximation falls below a certain tolerance level, all elements of the non-dominated set are guaranteed to be within a certain distance of the inner/outer approximation. For example, the \texttt{polySCIP}~\cite{Borndoerfer2016}, \texttt{Inner}~\cite{Csirmaz2021} and \texttt{Bensolve}~\cite{Lohne2017} softwares offer efficient implementations of Benson-type algorithms for multiobjective linear programming (MOLP). \texttt{Bensolve} also supports VO linear optimization. This website\footnote{\url{https://www.tu-ilmenau.de/en/mmor/software}} lists some available implementations in Python and/or Matlab of enclosure-based algorithms for nonlinear MO: e.g., AdEnA~\cite{Eichfelder2023a}, HyPaD~\cite{Eichfelder2024}, MOMIB~\cite{Eichfelder2024a} or MOMIX~\cite{Niebling2019} for non-convex MO (mixed-integer) nonlinear programming (MO(-MI)NLP). For more information about non-commercial available software solvers for MO/VO, the reader can refer to~\cite{Eichfelder2021b}.

Exploiting the structure of a problem is not free. Practitioners may need to dedicate effort to implementing their models in low-level code or the standard form required by a solver, which can limit adoption. Algebraic model languages such as CVX~\cite{Grant2008, Grant2014}, YALMIP~\cite{Lofberg2004}, AMPL~\cite{Fourer1990}, GAMS~\cite{Bussieck2004}, JuMP~\cite{Dunning2017} or Pyomo~\cite{Bynum2021Pyomo} remove this burden. They allow users to describe their problems in an algebraic and intuitive language that matches the original mathematical formulation. These languages handle the costly translation and pass the transformed optimization model to dedicated solvers. Few incorporate methods for MO/VO optimization problems. For example, GAMS (version \(50\)), a commercial mathematical modeling software, has a module for MO that provides a sandwiching algorithm and a rigid grid scalarization-based approach based on the augmented \(\varepsilon\)-constraint formulation for MOMILP and MOMINLP. \texttt{MultiObjectiveAlgorithms.jl}~\cite{Dowson2025}, based on JuMP, implements several scalarization-based algorithms for MOMILP and MOMINLP.

In the lineage of these works, this paper presents MOCVXPY, an open-source Python library for MO/VO convex optimization. MOCVXPY is built on top of CVXPY~\cite{Diamond2016}, a Python embedded-modeling language for single-objective convex optimization. Although other alternatives exist, e.g., Pyomo or PICOS~\cite{Sagnol2022}, several reasons explain this choice.
\begin{enumerate}
    \item CVXPY has an intuitive modeling API and is simple to use (just call the \texttt{solve()} method), though there is a slight decrease in computer time efficiency. MOCVXPY retains most of this API and strives to closely resemble the user interface of its single-objective counterpart library. Coded in Python, it can be incorporated into larger Python workflows.
    \item CVXPY can call several single-objective solvers and automatically assigns one according to the most specific subclass of the optimization problem provided. MOCVXPY takes advantage of this feature when solving its subproblems.
    \item CVXPY uses disciplined convex programming rules and reductions of atom functions to verify convexity of single-objective problems. MOCVXPY uses the same rules. This allows one to check that the MO/VO problem is convex (it may still be unbounded or infeasible). Otherwise, the library raises an exception and no MO/VO algorithm runs.
\end{enumerate}
On the contrary, MOCVXPY is not suited for more general categories of MO/VO problems: non-convex, mixed-integer. The reader whose problems belong to these classes can refer to the available implementations of AdeNa or HyPad, for example.

The remainder of this work is structured as follows. Section~\ref{section:Preliminaries} presents the problem and its associated main concepts. Section~\ref{section:A_toy_example} illustrates the features of the MOCVXPY library with a toy example. Section~\ref{section:The_algorithms} describes the three algorithms contained in the library and discusses some implementation details. 
This is followed by two real-world applications in Section~\ref{section:Real_world_applications} and concludes with a discussion in Section~\ref{section:Discussion}. 

\section{Preliminaries}\label{section:Preliminaries}

\subsection{Notations and conventions}

The organization of the different concepts described in this section is adapted from~\cite{Ararat2023, Ararat2022, Ararat2024, Dorfler2022}.

Throughout this work, let \(\mathbb{R}^q\) denote the \(q\)-dimensional Euclidean space with \(q \in \mathbb{N} \setminus \{0\}\), and \(\|.\|\) be the \(\ell_2\) norm on \(\mathbb{R}^q\). The \(\ell_2\) norm closed ball on \(\mathbb{R}^q\) centered at \(\bm{0} \in \mathbb{R}^p\) of radius \(\varepsilon > 0\) is denoted as \(\mathbb{B}(\bm{0}, \varepsilon) = \{\bm{z} \in \mathbb{R}^q: \|\bm{z}\| \leq \varepsilon \}\). Let \(\bm{1}_q\) denote the \(q\)-dimensional vector whose coordinates are equal to \(1\).

To compare vectors, the following notations are used. Given two vectors \(\bm{z}^1, \bm{z}^2 \in \mathbb{R}^q\)
\begin{itemize}
    \item \(\bm{z}^1 < \bm{z}^2 \Longleftrightarrow \bm{z}^1_i < \bm{z}^2_i\) for \(i = 1, 2, \ldots, q\);
    \item \(\bm{z}^1 \leq \bm{z}^2 \Longleftrightarrow \bm{z}^1_i \leq \bm{z}^2_i\) for \(i = 1, 2, \ldots, q\);
    \item \(\bm{z}^1 \lneq \bm{z}^2 \Longleftrightarrow \bm{z}^1 \leq \bm{z}^2\) and \(\bm{z}^1 \neq \bm{z}^2\).
\end{itemize}
Similarly, \(\bm{z}^1 > / \geq / \gneq \bm{z}^2 \Longleftrightarrow -\bm{z}^1 < / \leq / \lneq - \bm{z}^2\).

Let $\mathcal{A} \subseteq \mathbb{R}^q$. We denote by 
$\operatorname{cl}(\mathcal{A})$, $\operatorname{int}(\mathcal{A})$, 
$\operatorname{conv}(\mathcal{A})$, and $\operatorname{cone}(\mathcal{A})$ 
the closure, interior, convex hull, and conic hull of $\mathcal{A}$, respectively.
Recall that \(\operatorname{cone}(\mathcal{A}) = \{\lambda \bm{z}: \bm{z} \in \mathcal{A}, \lambda  \geq 0\}\). Given two non-empty sets \(\mathcal{A}, \mathcal{B} \subseteq \mathbb{R}^q\) and \(\lambda \in \mathbb{R}\), their Minkowski operations are defined by \(\mathcal{A} + \mathcal{B} = \{\bm{z}^1 + \bm{z}^2: \bm{z}^1 \in \mathcal{A}, \bm{z}^2 \in \mathcal{B}\}\), \(\lambda \mathcal{A} = \{\lambda \bm{z} : \bm{z} \in \mathcal{A}\}\) and \(\mathcal{A} - \mathcal{B} = \mathcal{A} + (-1) \mathcal{B}\). 

Let \(\mathcal{A} \subseteq \mathbb{R}^q\) be a non-empty set. For \(\bm{z} \in \mathbb{R}^q\), define \(d(\bm{z}, \mathcal{A}) = \displaystyle\inf_{\bm{z}' \in \mathcal{A}} \|\bm{z} - \bm{z}'\|\). Given a non-empty set \(\mathcal{B} \subseteq \mathbb{R}^q\), the Hausdorff distance between \(\mathcal{A}\) and \(\mathcal{B}\) is defined as
\[d_{H}(\mathcal{A}, \mathcal{B}) = \max \left\{\sup_{\bm{z}^1 \in \mathcal{A}} d(\bm{z}^1, \mathcal{B}), \sup_{\bm{z}^2 \in \mathcal{B}} d(\bm{z}^2, \mathcal{A})\right\}.\]

Recall that a convex polyhedral set \(\mathcal{A}\) can be expressed via an \(H\)-representation, i.e., as an intersection of a finite number of closed halfspaces
\begin{equation}
    \mathcal{A} = \bigcap_{j = 1}^l \{\bm{z} \in \mathbb{R}^q: (\bm{w}^j)^\top \bm{z} \geq \bm{\gamma}_j\} = \{\bm{z} \in \mathbb{R}^q: W \bm{z} \geq \bm{\gamma}\}
\end{equation}
where \(l \in \mathbb{N}\), \(W \in \mathbb{R}^{l \times q}\) whose rows are vectors \(\bm{w}^j \in \mathbb{R}^q\) for \(j = 1, \ldots, l\) and \(\bm{\gamma} \in \mathbb{R}^l\). Equivalently, \(\mathcal{A}\) can be also written via a \(V\)-representation, i.e.,
\begin{equation}
    \mathcal{A} = \operatorname{conv}\left(\{\bm{v}^1, \ldots, \bm{v}^s\}\right) + \operatorname{cone}\left(\{\bm{d^1}, \ldots, \bm{d}^r\}\right) = \{\bm{z} \in \mathbb{R}^q: \bm{z} = V \bm{\lambda} + D \bm{\mu}, \bm{1}_q^\top \bm{\lambda} = 1, \bm{\lambda}, \bm{\mu} \geq 0\}
\end{equation}
where \(V \in \mathbb{R}^{q \times s}\) whose columns \(\bm{v}^j \in \mathbb{R}^q\) for \(j = 1, \ldots, s\) are vertices of \(\mathcal{A}\), \(D \in \mathbb{R}^{q \times r}\) whose columns \(\bm{d}^j \in \mathbb{R}^q\) for \(j = 1, \ldots, r\) are extreme rays of \(\mathcal{A}\), \(\bm{\lambda} \in \mathbb{R}^s\) and \(\bm{\mu} \in \mathbb{R}^r\).

Let \(\mathcal{C} \subseteq \mathbb{R}^q\) be a convex cone. The \textit{dual cone} of \(\mathcal{C}\) is defined as \(\mathcal{C}^+ = \{\bm{w} \in \mathbb{R}^q: \forall \bm{z} \in \mathcal{C}, \bm{w}^\top \bm{z} \geq 0\}\): it is a closed convex cone. A cone is said to be \textit{solid} if \(\operatorname{int}(\mathcal{C}) \neq \emptyset\) and \textit{pointed} if \(-\mathcal{C} \cap \mathcal{C} = \{\bm{0}\}\). It is \textit{non-trivial} if \(\mathcal{C} \neq \emptyset\) and \(\mathcal{C} \neq \mathbb{R}^q\). If \(\mathcal{C} \subseteq \mathbb{R}^q\) is a solid pointed non-trivial convex cone, the relation \(\leq_\mathcal{C} = \{(\bm{z}^1, \bm{z}^2) \in \mathbb{R}^q: \bm{z}^2 - \bm{z}^1 \in \mathcal{C}\}\) defines a partial order on \(\mathbb{R}^q\)~\cite[Theorem 1.20]{Ehrgott2005MulticriteriaOptimization}. Given \(\bm{z}^1, \bm{z}^2 \in \mathbb{R}^q\), \(\bm{z}^1 \leq_{\mathcal{C}} \bm{z}^2\) is equivalent to \((\bm{z}^1, \bm{z}^2) \in \leq_{\mathcal{C}}\).

Let \(f: \mathbb{R}^n \rightarrow \mathbb{R}^q\). Given a closed convex solid and pointed cone \(\mathcal{C} \subseteq \mathbb{R}^q\), \(f\) is called \textit{\(\mathcal{C}\)-convex} if \(f\left(\lambda \bm{x}^1 + (1 - \lambda) \bm{x}^2\right) \leq_{\mathcal{C}} \lambda f(\bm{x}^1) + (1 - \lambda) f(\bm{x}^2)\) for all \(\bm{x}^1, \bm{x}^2 \in \mathbb{R}^n\) and \(\lambda \in [0, 1]\). Equivalently, \(f\) is \(\mathcal{C}\)-convex if \(\bm{x} \mapsto \bm{w}^\top f(\bm{x})\) is convex for all \(\bm{w} \in \mathcal{C}^+\). Given a set \(\mathcal{S} \subseteq \mathbb{R}^n\), \(f(\mathcal{S})\) denotes the set \(f(\mathcal{S}) = \{f(\bm{x}): \bm{x} \in \mathcal{S}\}\).

A set \(\mathcal{C} \subseteq \mathbb{R}^q\) is called a \textit{polyhedral convex cone} if there exists a matrix \(D \in \mathbb{R}^{q \times r}\) such that \(\mathcal{C} = \{D \bm{\mu}: \bm{\mu} \geq 0\} = \operatorname{cone} \left(\{\bm{d}^1, \ldots, \bm{d}^r\}\right)\) with \(\bm{d}^j \in \mathbb{R}^q\) for \(j = 1, \ldots, r\). Polyhedral convex cones have important properties~\cite{ Dorfler2022, Greer1984, Kaibel2011}:
\begin{enumerate}
    \item There exists a matrix \(W \in \mathbb{R}^{q \times l}\) such that \(\mathcal{C} = \{\bm{z} \in \mathbb{R}^q: W^\top \bm{z} \geq \bm{0}\}\). Specifically, \(\mathcal{C}^+ = \operatorname{cone}\left(\{\bm{w^1, \ldots, \bm{w}^l}\}\right)\) where \(\bm{w}^j\) is the vector that corresponds to the \(j\) column of \(W\).
    \item As for arbitrary closed convex cones, one has \((\mathcal{C}^+)^+ = \mathcal{C}\).
    \item Let \(W \in \mathbb{R}^{q\times l}\) such that \(\mathcal{C} = \{\bm{z} \in \mathbb{R}^q: W^\top \bm{z} \geq \bm{0}\}\). Then \(\mathcal{C}\) is pointed if and only if  \(\text{rank} \ W = q\).
    \item Suppose \(\mathcal{C}^+ \neq \{0\}\) is pointed. Let \(W \in \mathbb{R}^{q\times l}\) without zero columns, such that \(\mathcal{C} = \{\bm{z} \in \mathbb{R}^q: W^\top \bm{z} \geq \bm{0}\}\). Then \(\operatorname{int}(\mathcal{C}) = \{\bm{z} \in \mathbb{R}^q: W^\top \mathbf{z} > 0\}\).
\end{enumerate}

\subsection{Problem and optimal solutions}

MOCVXPY considers the \textit{convex vector optimization problem} (CVOP):
\begin{equation} \label{ref:VOP} \tag{VOP}
    \min_{\bm{x} \in \mathcal{X}} f(\bm{x}) = \left(f_1(\bm{x}), f_2(\bm{x}), \ldots, f_q(\bm{x})\right)^\top \text{ with respect to } \leq_{\mathcal{C}}
\end{equation}
where \(q \geq 2\) is the number of objectives of the problem. \(\mathcal{C} \subseteq \mathbb{R}^q\) is a pointed solid non-trivial polyhedral cone given as \(\mathcal{C} = \{\bm{z} \in \mathbb{R}^q: W^\top \bm{z} \geq 0\}\). The feasible set \(\mathcal{X} \neq \emptyset \subseteq \mathbb{R}^n\) is convex and \(f : \mathcal{X} \rightarrow \mathbb{R}^q\) is a \(\mathcal{C}\)-convex continuous function.
Note that for all \(\bm{z}^1, \bm{z}^2 \in \mathbb{R}^q\), \(\bm{z}^1 \leq_{\mathcal{C}} \bm{z}^2 \Longleftrightarrow W^\top \bm{z}^1 \leq W^\top \bm{z}^2\).

The upper image associated with \eqref{ref:VOP} plays a role in the characterization of a ``solution set''.

\begin{definition}[Adapted from~\cite{Dorfler2022}]
    The \textit{upper image} of \eqref{ref:VOP} is defined as
    \[\mathcal{P} = \operatorname{cl}(f(\mathcal{X}) + \mathcal{C}).\]
    Furthermore, \eqref{ref:VOP} is \textit{bounded} if there exists some \(\bm{z} \in \mathbb{R}^q\) such that \(\mathcal{P} \subseteq \{\bm{z}\} + \mathcal{C}\).
\end{definition}

\begin{figure}[!th]
    \centering
    \includegraphics{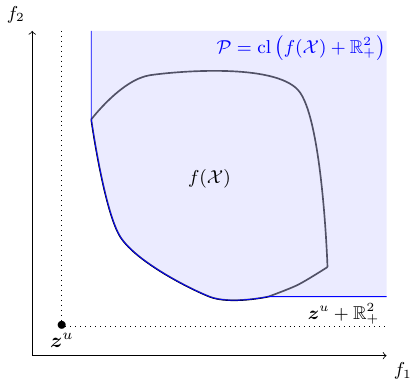}
    \caption{Illustration of the upper image of a CVOP with \(\mathcal{C} = \mathbb{R}^2_+\).}
    \label{fig:upper_image}
\end{figure}

Clearly, \(\mathcal{P}\) is a closed convex set. Figure~\ref{fig:upper_image} shows an upper image of a VOP instance with order cone $\mathcal{C} = \mathbb{R}^2_+$.
Note that this problem is bounded since \(\mathcal{P} \subset \bm{z}^u + \mathbb{R}^2_+\). 

This work uses the following assumption.
\begin{assumption}
    \eqref{ref:VOP} is bounded.
\end{assumption}
Note that, if the feasible set \(\mathcal{X}\) of \eqref{ref:VOP} is compact, then by continuity of \(f\), \(f(\mathcal{X})\) is compact as well and the problem is already bounded~\cite[Theorem~3.2]{Dorfler2022a}.

The problem associated with \(\mathcal{C} = \mathbb{R}_+^q\) is called a \textit{convex multiobjective optimization problem}, and is written as:
\begin{equation} \label{ref:MOP} \tag{MOP}
    \min_{\bm{x} \in \mathcal{X}} f(\bm{x}) = \left(f_1(\bm{x}), f_2(\bm{x}), \ldots, f_q(\bm{x})\right)^\top.
\end{equation}
In this case, the corresponding ordering relation \(\leq_{\mathbb{R}^q_+}\), also called the \textit{natural ordering}, is omitted. This special case, denoted \eqref{ref:MOP}, is important enough that we spend some time on it. The following optimality notions are defined for \eqref{ref:MOP}.

\begin{definition}[Adapted from~\cite{Eichfelder2021b}]
Let \(\mathcal{C} = \mathbb{R}^q_+\). Let \(\varepsilon > 0\).
\begin{itemize}
    \item A decision vector \(\bar{\bm{x}} \in \mathcal{X}\) is called a \textit{Pareto-efficient solution/point} for \eqref{ref:MOP}, and its image \(f(\bar{\bm{x}})\) a \textit{non-dominated point} for \eqref{ref:MOP}, if there is no \(\bm{x} \in \mathcal{X}\) such that \(f(\bm{x}) \lneq f(\bar{\bm{x}})\).
    \item A decision vector \(\bar{\bm{x}} \in \mathcal{X}\) is called a \textit{weakly Pareto-efficient solution/point} for \eqref{ref:MOP}, and its image \(f(\bar{\bm{x}})\) a \textit{weakly non-dominated point} for \eqref{ref:MOP}, if there is no \(\bm{x} \in \mathcal{X}\) such that \(f(\bm{x}) < f(\bar{\bm{x}})\).
    \item A decision vector \(\bar{\bm{x}} \in \mathcal{X}\) is called a \textit{\(\varepsilon\)-Pareto-efficient solution/point} for \eqref{ref:MOP}, and its image \(f(\bar{\bm{x}})\) a \(\varepsilon\)-non-dominated point for \eqref{ref:MOP}, if there is no \(\bm{x} \in \mathcal{X}\) such that \(f(\bm{x}) \lneq f(\bar{\bm{x}}) - \varepsilon\bm{1}_q\).
\end{itemize}
\end{definition}

Any Pareto-efficient solution for \eqref{ref:MOP} is also weakly efficient for \eqref{ref:MOP}. In practice, Pareto-efficient solutions are more desirable than weakly efficient ones because, for the latter, it is possible to decrease one objective component without deteriorating the others~\cite{Eichfelder2021b}. However, weakly efficient solutions play a role in  the mathematical characterization of outputs generated by many numerical algorithms for \eqref{ref:MOP}.

The set of all Pareto-efficient solutions is called the \textit{Pareto(-efficient) set} for \eqref{ref:MOP} and its associated image set by the objective function is the \textit{non-dominated set} for \eqref{ref:MOP} or \textit{Pareto front}, denoted as \(\mathcal{Z}_P\). The set of \(\varepsilon\)-non dominated points for \eqref{ref:MOP} is denoted as \(\mathcal{Z}^{\varepsilon}_P\). Note that \(\mathcal{Z}_P \subseteq \mathcal{Z}^{\varepsilon}_P\).

\begin{figure}[!th]
    \centering
    \includegraphics[width=0.5\linewidth]{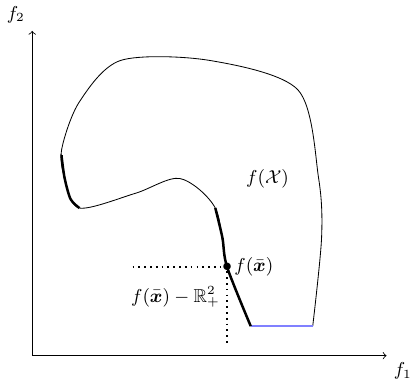}
    \caption{Illustration of the set of (weakly) Pareto efficient solutions for a biobjective non-convex MOP. The weakly non-dominated set is the union of the non-dominated set indicated in very thick black, and another part of the objective space, drawn in blue.}
    \label{fig:weakly_and_efficient_solutions}
\end{figure}

In terms of set notation, one can write that \(\bar{\bm{x}} \in \mathcal{X}\) is Pareto-efficient for \eqref{ref:MOP} if and only if
\[\left(\{f(\bar{\bm{x}})\} - \mathbb{R}^q_+\right) \cap f(\mathcal{X}) = \{f(\bar{\bm{x}})\}.\]
Similarly, \(\bar{\bm{x}} \in \mathcal{X}\) is weakly Pareto-efficient for \eqref{ref:MOP} if and only if
\[\left(\{f(\bar{\bm{x}})\} - \operatorname{int}\left(\mathbb{R}^q_+\right)\right) \cap f(\mathcal{X}) = \emptyset.\]
Figure~\ref{fig:weakly_and_efficient_solutions} shows the non-dominated set for a non-convex biobjective MOP. The image of the weakly efficient solutions of this problem is the union of the non-dominated set, in very thick lines and another part of the feasible objective space, shown in blue.

With this set notation, the concept of a (weakly) efficient solution can be generalized to vector optimization. 
This is possible because the ordering cone \(\mathcal{C}\) is polyhedral, solid, and pointed.

\begin{definition}
    \begin{itemize}
    \item A decision vector \(\bar{\bm{x}} \in \mathcal{X}\) is called a \textit{minimizer} for \eqref{ref:VOP}, and its image \(f(\bar{\bm{x}})\) a \textit{\(\mathcal{C}\)-minimal point} of \(f(\mathcal{X})\), if
    \[\left(\{f(\bar{\bm{x}})\} - \mathcal{C} \right) \cap f(\mathcal{X})= \{f(\bar{\bm{x}})\}.\]
    \item A decision vector \(\bar{\bm{x}} \in \mathcal{X}\) is called a \textit{weak minimizer} for \eqref{ref:VOP}, and its image \(f(\bar{\bm{x}})\) a \textit{weakly \(\mathcal{C}\)-minimal point} of \(f(\mathcal{X})\), if
    \[\left(\{f(\bar{\bm{x}})\} - \operatorname{int} \left(\mathcal{C}\right) \right) \cap f(\mathcal{X})= \emptyset.\]
    \end{itemize}
\end{definition}

With this terminology, a (weakly) non-dominated point for (MOP) corresponds to a (weakly) \(\mathbb{R}^q_+\)-minimal point of \(f(\mathcal{X})\).

Given a parameter \(\bm{w} \in \mathcal{C}^+\), the \textit{weighted-sum scalarization} of \eqref{ref:VOP} is defined as
\begin{equation} \label{ref:WS} \tag{WS(\(\bm{w}\))}
    \min_{\bm{x} \in \mathcal{X}} \bm{w}^\top f(\bm{x}).
\end{equation}

The following well-known proposition recalls the connection between (\ref{ref:WS}) and weak minimizers of \eqref{ref:VOP}.

\begin{proposition}{(\cite[Corollary 5.29]
{Jahn2009VectorOptimization}).}\label{prop:weighted_sum_correspondance} Let \(\bm{w} \in \mathcal{C}^+ \setminus \{\bm{0}\}\). Then, every optimal solution of (\ref{ref:WS}) is a weak minimizer of \eqref{ref:VOP}. Conversely, for each weak minimizer \(\bm{x} \in \mathcal{X}\) of \eqref{ref:VOP}, there exists \(\bm{w} \in \mathcal{C}^+ \setminus \{\bm{0}\}\) such that \(\bm{x}\) is an optimal decision of (\ref{ref:WS}).
\end{proposition}

In multiobjective as in vector optimization, the user is generally not satisfied with obtaining only one efficient solution/minimizer. Rather, he looks for a set of different feasible decision vectors whose images by the objective function represent the best trade-offs according to the ordering cone \(\mathcal{C}\). The characterization of such a set follows the one presented in~\cite{Lohne2011Vector}.

\begin{definition}[Adapted from~\cite{Heyde2011}]\label{def:infimizer_VOP}
A non-empty set \(\bar{\mathcal{X}} \subseteq \mathcal{X}\) is called an \textit{infimizer} of \eqref{ref:VOP} if \(\operatorname{cl}\left( \operatorname{conv}\left(f(\bar{\mathcal{X}}) + \mathcal{C}\right) \right) = \mathcal{P}\). An infimizer \(\bar{\mathcal{X}}\) of \eqref{ref:VOP} is called a \textit{(weak) solution set} of \eqref{ref:VOP} if it consists of only (weak) minimizers.
\end{definition}

\begin{figure}[!th]
    \centering
    \includegraphics[]{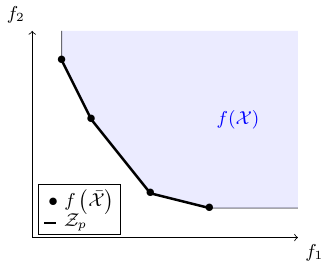}
    \caption{An illustration of a minimal infimizer for a biobjective minimization CVOP. Here, \(\mathcal{P} = f(\mathcal{X})\).}
    \label{fig:infimizer_example}
\end{figure}

\begin{remark}
    When \(\mathcal{C} = \mathbb{R}^q_+\), a solution set in the sense of Definition~\ref{def:infimizer_VOP} is a subset of the Pareto set, but is not necessarily equal to it, as illustrated in Figure~\ref{fig:infimizer_example}. Thus, a solution set can be interpreted as a ``minimum'' number of solutions in the feasible decision space whose images ``capture the shape of the Pareto front''.
\end{remark}

For the reasons discussed above, it can be difficult or even impossible to compute a solution set in the sense of Definition~\ref{def:infimizer_VOP}. One looks for a finite approximate solution set of (\ref{ref:VOP}, that can represent the true solution set.

\begin{definition}[Adapted from~\cite{Dorfler2022}]
    Suppose that \eqref{ref:VOP} is bounded. Let \(\varepsilon > 0\). Let \(\bar{\mathcal{X}} \subseteq \mathcal{X}\) be a non-empty finite set and define \(\bar{\mathcal{P}} = \operatorname{conv}(f(\bar{\mathcal{X}})) + \mathcal{C}\). The set \(\mathcal{\bar{X}}\) is called a \textit{finite \(\varepsilon\)-infimizer} of \eqref{ref:VOP} if \(\bar{\mathcal{P}} + \mathbb{B}(\bm{0}, \varepsilon) \supseteq \mathcal{P}\). The set \(\mathcal{\bar{X}}\) is called a \textit{finite (weak) \(\varepsilon\)-solution set} of \eqref{ref:VOP} if it is an \(\varepsilon\)-infimizer that consists of only (weak) minimizers.
\end{definition}

\begin{figure}[!th]
    \centering
    \includegraphics{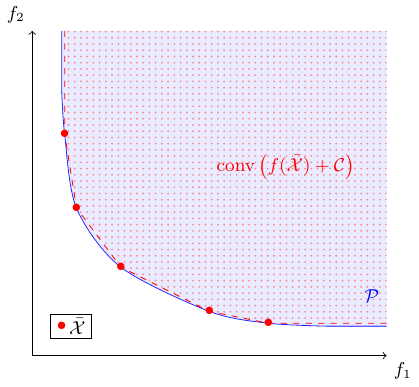}
    \caption{Illustration of an \(\varepsilon\)-infimizer, inspired by~\cite{Dorfler2022}. Here, \(\mathcal{C} = \mathbb{R}^2_+\).}
    \label{fig:epsilon_infimizer}
\end{figure}

An illustration of an \(\varepsilon\)-infimizer is given in Figure~\ref{fig:epsilon_infimizer} for a biobjective CVOP.

Note that \(\bar{\mathcal{P}} + \mathbb{B}(\bm{0}, \varepsilon) \supseteq \mathcal{P} \supseteq \bar{\mathcal{P}}\). Thus, a finite \(\varepsilon\)-infimizer (or a finite \(\varepsilon\)-solution set) \(\bar{\mathcal{X}}\) defines outer and inner approximations of the upper image \(\mathcal{P}\). 
The Hausdorff distance can be measured between these two sets.

In multiobjective optimization, one can exploit the natural ordering symmetry to ``sandwich'' the non-dominated set of \eqref{ref:MOP} between a suitable lower bound set and a suitable upper bound set. The combination of these approximation is called an enclosure for the non-dominated set \(\mathcal{Z}_P\)~\cite{Eichfelder2021a, Eichfelder2023a}.

\begin{definition}[Adapted from~\cite{Eichfelder2023a}] Let \(\mathcal{C} = \mathbb{R}^q_+\).
Let \(L, U \subset \mathbb{R}^q\) be two non-empty finite sets with \(\mathcal{Z}_P \subseteq L + \mathbb{R}^q_+\) and \(\mathcal{Z}_P \subseteq U- \mathbb{R}^q_+\).
Then \(L\) is called a \textit{lower bound set}, \(U\) is called an \textit{upper bound set}, and the set \(\mathcal{A}(L, U)\) defined as
\[\mathcal{A}(L, U) = (L + \mathbb{R}^q_+) \cap (U - \mathbb{R}^q_+) = \bigcup_{\begin{subarray}{c} (\bm{l}, \bm{u}) \in L \times U, \\ \bm{l} \leq \bm{u} \end{subarray}} [\bm{l}, \bm{u}]\]
with \([\bm{l}, \bm{u}] = \{\bm{z} \in \mathbb{R}^q: \bm{l} \leq \bm{z} \leq \bm{u}\}\) is called an \textit{enclosure} (or a box approximation) of the non-dominated set \(\mathcal{Z}_P\) of \eqref{ref:MOP} given \(L\) and \(U\).
\end{definition}

Note that the set \(L + \mathbb{R}^q_+\) is an outer approximation of \(\mathcal{P}\) for \eqref{ref:MOP} and \(U + \mathbb{R}^q_+\) is an inner approximation \(\mathcal{P}\) for \eqref{ref:MOP} if \(U \subset \mathcal{P}\).

One can measure the quality of an enclosure \(\mathcal{A}(L, U)\) by its \textit{width} \(w\left(\mathcal{A}(L, U)\right)\), introduced in~\cite{Eichfelder2021a}. It is defined by the optimal value
\[w\left(\mathcal{A}(L, U)\right) = \max_{\bm{l} \in L, \bm{u} \in U} \{s(\bm{l}, \bm{u}): \bm{l} \leq \bm{u}\}\]
where \(s(\bm{l}, \bm{u}) = \displaystyle\min_{1\leq i \leq q} u_i - l_i\) is the shortest edge of the box \([\bm{l}, \bm{u}]\). The width of an enclosure is often used as a stopping criterion for enclosure-based approaches, e.g.,~\cite{Eichfelder2021a, Eichfelder2021, Eichfelder2023a}. A detailed discussion of this quality measure can be found in~\cite{Eichfelder2021b,Eichfelder2021a, Eichfelder2021}. In particular, the following property~\cite[Lemma 3.1]{Eichfelder2021a} justifies the definition of this quality measure: if $w\left(\mathcal{A}(L, U)\right) < \varepsilon$ for some \(\varepsilon > 0\), then:
\[\mathcal{A}(L, U) \cap f(\mathcal{X}) \subseteq \mathcal{Z}_P^\varepsilon.\]

The computation and the updating of an enclosure during the optimization of \eqref{ref:MOP} strongly depend on the choice of the lower bound set \(L\) and the upper bound set \(U\). A common strategy for computing an upper bound set is to use so-called local upper bound sets, which were initially introduced in~\cite{Klamroth2015} and generalized in~\cite{Eichfelder2023a}. These sets are related to search regions in the objective space, where new solutions could be found. Before presenting the definition, a set of objective vectors \(\emptyset \neq N \subseteq \mathbb{R}^q\) is said to be \textit{stable} if for any \(\bm{z}^1, \bm{z}^2 \in N\), either \(\bm{z}^1 = \bm{z}^2\) or \(\bm{z}^1 \nleq \bm{z}^2\).

\begin{definition}[Adapted from \cite{Eichfelder2023a}] Let \(\mathcal{B} \subseteq \mathbb{R}^q\) be a set, referred to as a zone of interest in the objective space, with \(\text{int}(\mathcal{B}) \neq \emptyset\). Let \(N \subset \mathbb{R}^q \neq \emptyset\) be a finite and stable set. Then the lower search region for \(N\) is \(s(N) = \{\bm{z} \in \text{int}(\mathcal{B}): \bm{z}' \nleq \bm{z} \text{ for every } \bm{z}' \in N\}\) and the lower search zone for some \(\bm{u} \in \mathbb{R}^q\) is \(c(\bm{u}) = \{\bm{z} \in \text{int}(\mathcal{B}): \bm{z} < \bm{u}\}\). 
A set \(U = U(N) \subseteq \mathbb{R}^q\) is called a \textit{local upper bound set} with respect to \(N\)
if
\begin{enumerate}
    \item \(s(N) = \bigcup_{\bm{u} \in U(N)} c(\bm{u})\),
    \item \(\bm{u}^1 - \text{int}(\mathbb{R}^q_+) \nsubseteq \bm{u}^2 - \text{int}(\mathbb{R}^q_+)\) for all \(\bm{u}^1, \bm{u}^2 \in U(N)\) with \(\bm{u}^1 \neq \bm{u}^2\).
\end{enumerate}
Each point \(\bm{u} \in U(N)\) is called a \textit{local upper bound}.    
\end{definition}

Similarly, one can define a local lower bound set for the stable set \(N\).

\begin{definition}[Adapted from \cite{Eichfelder2023a}] Let \(\mathcal{B} \subseteq \mathbb{R}^q\) be a set, referred to as a zone of interest in the objective space, with \(\operatorname{int}(\mathcal{B}) \neq \emptyset\). Let \(N \subset \mathbb{R}^q \neq \emptyset\) be a finite and stable set. Then the upper search region for \(N\) is \(s(N) = \{\bm{z} \in \text{int}(\mathcal{B}): \bm{z}' \ngeq \bm{z} \text{ for every } \bm{z}' \in N\}\) and the upper search zone for some \(\bm{l} \in \mathbb{R}^q\) is \(c(\bm{l}) = \{\bm{z} \in \text{int}(\mathcal{B}): \bm{z} > \bm{l}\}\). A set \(L = L(N) \subseteq \mathbb{R}^q\) is called a \textit{local lower bound set} with respect to \(N\) if
\begin{enumerate}
    \item \(s(N) = \bigcup_{\bm{l} \in L(N)} c(\bm{l})\),
    \item \(\bm{l}^1 + \text{int}(\mathbb{R}^q_+) \nsubseteq \bm{l}^2 + \text{int}(\mathbb{R}^q_+)\) for all \(\bm{l}^1, \bm{l}^2 \in L(N)\) with \(\bm{l}^1 \neq \bm{l}^2\).
\end{enumerate}
Each point \(\bm{l} \in L(N)\) is called a \textit{local lower bound}.    
\end{definition}

\begin{figure}[!th]
    \centering
    \includegraphics{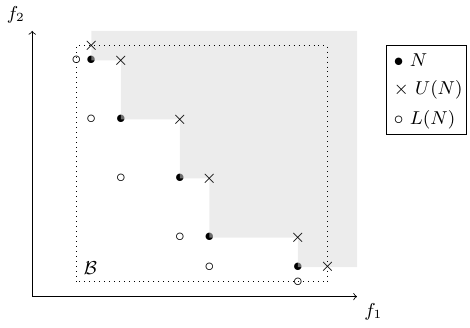}
    \caption{Illustration of local lower and upper bound sets for a biobjective minimization problem. The gray zone refers to the zone in the objective space dominated by the stable set \(N\).}
    \label{fig:llb_and_lub}
\end{figure}

Figure~\ref{fig:llb_and_lub} illustrates the notion of local lower and upper bound sets for a biobjective minimization problem. Interested readers can refer to~\cite{Dachert2017, Eichfelder2023a, Klamroth2015} for existing procedures to compute such sets. Under certain assumptions, the local lower bound and upper bound sets given a finite stable set \(N\) are guaranteed to be finite.

\begin{remark}
    In~\cite{Eichfelder2025}, the authors extend the concept of a local upper bound set to polyhedral ordering cones. However, they show that such a set can be infinite, hence difficult to use in numerical algorithms.
\end{remark}

\section{A toy example}\label{section:A_toy_example}

In this section, this biobjective minimization problem illustrates the functioning of MOCVXPY:
\begin{equation}\label{ex:disc}
\begin{array}{ll}
\displaystyle\min_{\bm{x} \in \mathbb{R}^2} & f(\bm{x}) = (x_1, x_2)^\top \\
\text{subject to} & (x_1 - 1)^2 + (x_2 - 1)^2 \leq 1 \\
& \bm{x} \geq 0.
\end{array} 
\end{equation}
The following code constructs the problem.
\begin{minted}[frame=lines, linenos]{python}
import cvxpy as cp
import mocvxpy as mocp

n = 2
x = mocp.Variable(n)
objectives = [cp.Minimize(x[0]), cp.Minimize(x[1])]
constraints = [cp.sum_squares(x - 1) <= 1,
               x >= 0]
pb = mocp.Problem(objectives, constraints)
\end{minted}

MOCVXPY reuses most of CVXPY's syntax. Thus, users can describe their convex problems using the same atomic functions implemented in CVXPY. The two main differences with CVXPY are the use of a \texttt{mocvxpy.Variable} instance instead of a \texttt{cvxpy.Variable} instance; and a \texttt{mocvxpy.Problem} instance instead of a \texttt{cvxpy.Problem} instance. The \texttt{mocvxpy.Variable} class is built on the \texttt{cvxpy.Variable} class and contains additional fields that will store optimal decision values once the \texttt{mocvxpy.Problem} is solved.

The following code launches the optimization problem.
\begin{minted}[frame=lines]{python}
objective_values = pb.solve()
\end{minted}
If the problem is solved, the set of optimal objective vectors obtained during optimization is stored in \texttt{objective\_values}. This follows the API of CVXPY, in which the method \texttt{solve()} of \texttt{cvxpy.Problem} returns the optimal value of a single-objective problem. The \texttt{mocvxpy.Problem} also populates some fields containing additional information about the optimal solutions. For example, all optimal values for the objectives, constraints and decision variables are accessible via the \texttt{.values} field.
\begin{minted}[frame=lines, linenos]{python}
print(x.values) # Display all optimal decision vectors found by an algorithm.
# Access to the second optimal decision vector is done with
print(x.values[1])
# Display its associated constraint values and objective values.
print(constraints[0].values[1], constraints[1].values[1])
print(objectives[0].values[1], objectives[1].values[1]) 
\end{minted}
All solutions obtained during the optimization are weak minimizers for \eqref{ref:VOP} or weakly Pareto-efficient for \eqref{ref:MOP}. By Proposition~\ref{prop:weighted_sum_correspondance}, each solution is an optimal solution of the weighted-sum scalarization (\ref{ref:WS}) for some weight values \(\bm{w} \in \mathcal{C}^+\). These weights are stored in the \texttt{.dual\_values} field of the objectives list after solving the problem. 
\begin{minted}[frame=lines, linenos, mathescape=true]{python}
# Display all weight values associated to optimal solutions of the problem.
print(objectives[0].dual_values, # $w_1$
      objectives[1].dual_values) # $w_2$
\end{minted}
Similarly, MOCVXPY provides two functions that compute the local lower and upper bound sets of a finite set of optimal objective vectors.
\begin{minted}[frame=lines, linenos, mathescape=true]{python}
# Display the local lower bound set associated to objective_values.
print(mocp.local_lower_bounds(objective_values))
# Display the local upper bound set associated to objective_values.
print(mocp.local_upper_bounds(objective_values))
\end{minted}

\begin{figure}[!th]
    \centering
    \includegraphics[scale=0.6]{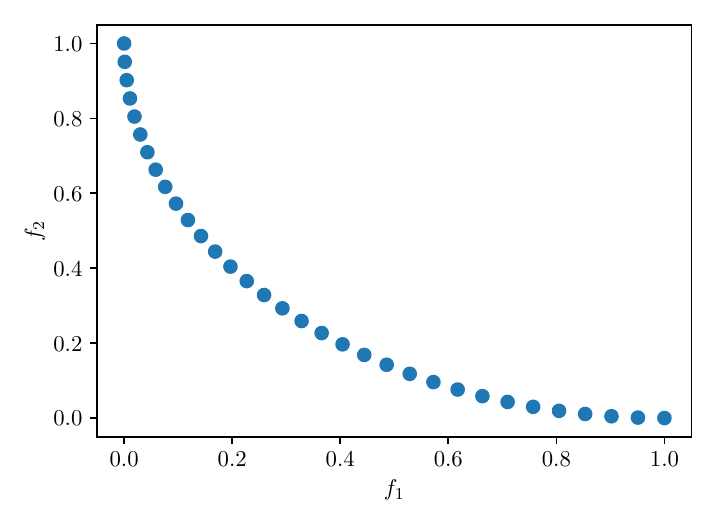}
    \caption{Approximation of the non-dominated set of Problem~\ref{ex:disc}, computed with default options.}
    \label{fig:toy_problem}
\end{figure}

MOCVXPY seamlessly integrates with other Python packages. For instance, although MOCVXPY lacks plotting utilities, one can use Matplotlib to display the optimal objective vectors.
\begin{minted}[frame=lines, linenos]{python}
from matplotlib import pyplot as plt

ax = plt.figure().add_subplot()
ax.scatter(
    [f1_value for f1_value in objective_values[:,0]],
    [f2_value for f2_value in objective_values[:,1]],
)
ax.set_xlabel("$f_1$")
ax.set_ylabel("$f_2$")
plt.show()
\end{minted}
The Pareto front approximation is shown on Figure~\ref{fig:toy_problem}.

To solve a convex vector optimization problem with MOCVXPY, the user must define a non-trivial, pointed, solid polyhedral ordering cone. 
MOCVXPY accepts polyhedral cones described by their H-representation or their V-representation (in the latter case, it is given as a conic combination of extreme rays). If no cone is specified, MOCVXPY assumes the problem is multiobjective and uses the non-negative orthant \(\mathbb{R}^q_+\) as the ordering cone.

\begin{figure}[!th]
    \centering
    \includegraphics[scale=0.6]{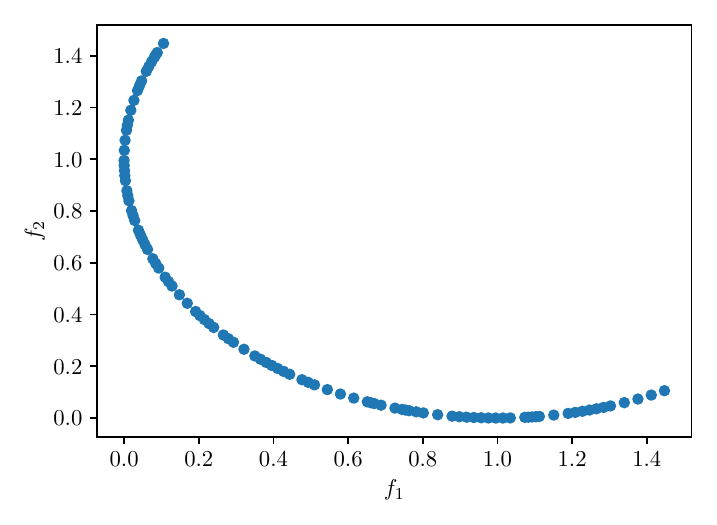}
    \caption{Approximation of the objective values of a solution of Problem~\ref{ex:disc} with order cone \(\mathcal{C}_1 = \text{cone}\left(\left\{(1, 2)^\top, (2, 1)^\top\right\}\right)\), computed with default options.}
    \label{fig:toy_problem_cone}
\end{figure}

The following code snippet details how to model and solve Problem~\ref{ex:disc} under the partial ordering relation \(\leq_{\mathcal{C}_1}\) taken from~\cite{Ararat2022}, where \(\mathcal{C}_1\) is defined as
\[\mathcal{C}_1 = \text{cone}\left(\left\{(1, 2)^\top, (2, 1)^\top\right\}\right).\]
\begin{minted}[frame=lines, linenos]{python}
import cvxpy as cp
import mocvxpy as mocp
import numpy as np

n = 2
x = mocp.Variable(n)
objectives = [cp.Minimize(x[0]), cp.Minimize(x[1])]
constraints = [cp.sum_squares(x - 1) <= 1,
               x >= 0]
C1 = mocp.compute_order_cone_from_its_rays(np.array([[1, 2], [2, 1]]))
pb = mocp.Problem(objectives, constraints, C1)

# Solve the problem.
objective_values = pb.solve()
\end{minted}
On line 10, the extreme rays of the cone are passed as a two-dimensional NumPy array (each row of the matrix represents a ray) to the \texttt{compute\_order\_cone\_from\_its\_rays} function. Behind the scenes, the function creates a cone described by its H-representation, which the algorithms will use during the resolution. On line 11, a \texttt{mocxvpy.Problem} instance is created using the \texttt{constraints} and \texttt{objectives} lists and the \texttt{C1} cone. Except for these two differences, modeling and solving the problem is the same as for its previous multiobjective variant. Figure~\ref{fig:toy_problem_cone} shows the obtained approximation of the objective values of a solution of Problem~\ref{ex:disc}.

\section{The algorithms}\label{section:The_algorithms}

\subsection{General presentation}

Table~\ref{tab:algorithms_summary} summarizes the algorithms supported by MOCVXPY. The library currently supports three algorithms. The Multiobjective Optimization by Norm Min Optimization (MONMO) algorithm is an outer approximation algorithm based on a norm min scalarization~\cite{Ararat2022}. The Multiobjective Optimization by Vertex Selection (MOVS) algorithm~\cite{Dorfler2022} alternates between two steps at each iteration: first, a vertex is selected from the current polyhedral outer approximation; then, a Pascoletti-Serafini scalarization subproblem is solved. These two solvers can tackle problems of the form \eqref{ref:VOP}. AdEnA~\cite{Eichfelder2023a} is an enclosure-based framework that refines an enclosure of the non-dominated set along the iterations using a Pascoletti-Serafini scalarization step.
\begin{table}[!ht]
    \centering
    \begin{tabular}{llccll}
    \hline
    Algorithm & Article & VO & Parallelism & Category & Solution\\
    \hline \hline
    AdEnA & \cite{Eichfelder2023a} & & & Enclosure-based & finite weakly efficient set\\
    \hline
    MONMO & \cite{Ararat2023, Ararat2022, Ararat2024} & \cmark & \cmark & Benson-type & \(\varepsilon\)-solution set\\
    \hline
    MOVS & \cite{Dorfler2022} & \cmark & \cmark &  Benson-type & \(\varepsilon\)-solution set \\
    \hline
    \end{tabular}
    \caption{Multiobjective optimization algorithms supported by MOCVXPY.}
    \label{tab:algorithms_summary}
\end{table}

The user can specify which algorithm to use as an argument of the \texttt{.solve()} method. The following code gives an example of how to use the MONMO solver.
\begin{minted}[frame=lines]{python}
pb.solve(solver="MONMO")
\end{minted}
MOVS is the default solver.

\begin{figure}[!th]
    \centering
    \subfigure[AdEnA]{
    \includegraphics[scale=0.5]{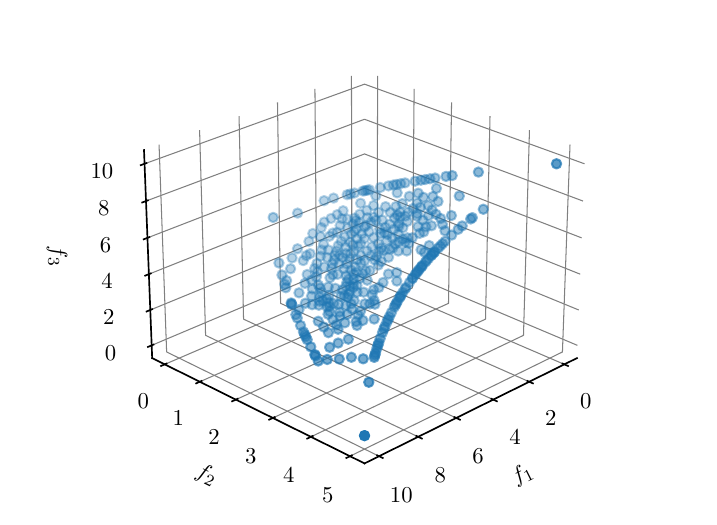}
    }
    \quad
    \subfigure[MONMO]{
    \includegraphics[scale=0.5]{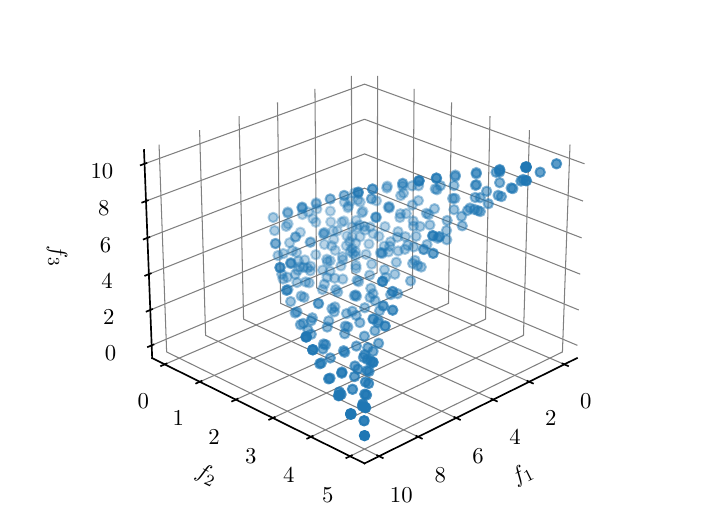}
    }
    \quad
    \subfigure[MOVS]{
    \includegraphics[scale=0.5]{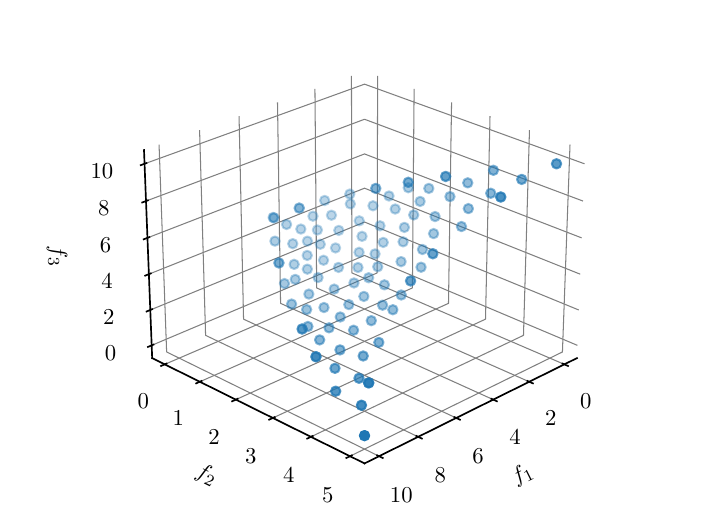}
    }
    \caption{Pareto front approximations using AdEnA, MONMO and MOVS on the following three-objective problem, taken from~\cite{Ararat2022}:
    \(\displaystyle\min_{\bm{x} \in \mathbb{R}^2} f(\bm{x}) = \left( \|\bm{x} - \bm{a}^1\|_2^2, \|\bm{x} - \bm{a}^2\|_2^2, \|\bm{x} - \bm{a}^3\|_2^2\right)^\top\) with respect to \(\leq_{\mathbb{R}^3_+}\), subject to \(x_1 + x_2 \leq 10, x_1 \in [0, 10], x_2 \in [0, 4]\); where \(\bm{a}^1 = (1, 1)^\top\), \(\bm{a}^2 = (2, 3)^\top\) and \(\bm{a}^3 = (4, 2)^\top\).}
    \label{fig:three_objective_mo_pb}
\end{figure}

Empirically, MOVS produces smaller solution sets than MONMO (related to the number of scalarization subproblems solved), whose images by the objective function are more evenly distributed in the objective space. In terms of computation time, MOVS should outperform MONMO for large \eqref{ref:MOP}/\eqref{ref:VOP} instances. Indeed, in this case, MOVS will spend most of its time part solving the scalarization subproblems and less time on the vertex selection procedure, which does not exist for MONMO. The AdEnA implementation in MOCVXPY is also computationally efficient. One must keep in mind that the primary goal of this algorithm is to compute a precise enclosure of \(\mathcal{Z}_p\) for \eqref{ref:MOP}; the solution set is a byproduct. Figure~\ref{fig:three_objective_mo_pb} shows the solution sets obtained by the different algorithms for a three-objective problem. We encourage practitioners to test the different algorithms according to their application.

\begin{remark}
    The open-source AdEnA implementation, available at \url{https://github.com/LeoWarnow/AdEnA}, is based on the Pyomo optimization modeling language. Describing a mathematical model with this implementation requires more work than with MOCVXPY, but it is more powerful because it can handle non-convex and/or mixed-integer multiobjective problems. Conversely, MOCVXPY benefits from the CVXPY environment which primarily targets disciplined convex problems. When the application belongs to this domain, faster computational time may be observed, particularly when passing the model of a scalarization subproblem to a single-objective convex solver.
\end{remark}

All of these algorithms have convergence guarantees. They are also parameter-free in the sense that they do not require the user to select an external mathematical parameter for the scalarization subproblems, e.g., a direction parameter as in~\cite{Keskin2023, Lohne2014}. The parameterized subproblems involved in these methods are, by construction, invariant with respect to the CVXPY framework. For a feasible, continuous and bounded \eqref{ref:VOP}/\eqref{ref:MOP} expressed in CVXPY's algebraic language, the subproblems are convex and can be handled by a conic solver. If solved to optimality, they satisfy strong duality. This property is not satisfied, for example in~\cite{Raimundo2020}, where a mixed-integer linear subproblem is employed.

\subsection{Scalarization formulations}\label{sec:scalarization_formulations}

Table~\ref{tab:scalarization_formulations} lists the scalarization formulations used by the different algorithms.

\begin{table}[!ht]
    \centering
    \begin{tabular}{p{1.5cm}p{5cm}p{3cm}p{3.5cm}}
    \hline
    Name & Formulation & Parameters & Algorithms \\
    \hline \hline
    (\ref{ref:WS}) & \begin{math}\displaystyle\min_{\bm{x} \in \mathcal{X}} \bm{w}^\top f(\bm{x})\end{math} & weights \(\bm{w} \in \mathbb{R}^q\) & AdeNa, MOVS, MONMO (initialization step) \\
    \hline
    (PS(\(\bm{v}, \bm{d}\))) & \begin{math} \begin{array}{ll} \min & t \\ 
    \text{subject to} & f(\bm{x}) \leq_{\mathcal{C}} \bm{v} + t \bm{d}\\ & \bm{x} \in \mathcal{X}, t \in \mathbb{R}\end{array} \end{math} & target \(\bm{v} \in \mathbb{R}^q\), direction \(\bm{d} \in \mathbb{R}^q\) & AdeNa, MOVS \\
    \hline
    (NM(\(\bm{v}\))) & \begin{math} \begin{array}{ll} \min & \| \bm{z} \| \\ 
    \text{subject to} & f(\bm{x}) \leq_{\mathcal{C}} \bm{z} + \bm{v}\\ & \bm{x} \in \mathcal{X}, \bm{z} \in \mathbb{R}^q \end{array} \end{math} & target \(\bm{v} \in \mathbb{R}^q\) & MONMO \\
    \hline
    \end{tabular}
    \caption{Scalarization formulations involved in MOCVXPY.}
    \label{tab:scalarization_formulations}
\end{table}

Although it is not necessary for practitioners to understand how the different algorithms use these scalarization formulations (the interested reader can refer to the references listed in Table~\ref{tab:algorithms_summary} for a complete mathematical description of these algorithms), MOCVXPY allows users to customize solver options for subproblem resolution. Thus, knowing the scalarization formulation associated with a given solver enables the user to consider the structure of its problem when choosing a solver. Otherwise, the practitioner can let CVXPY handle this choice. All single-objective solvers supported by CVXPY, along with their respective options \footnote{They can be found at \url{https://www.cvxpy.org/tutorial/solvers/index.html}.}, can be passed to MOCVXPY algorithms. However, this does not guarantee that these solvers will always solve the subproblems. The MO/VO methods also have their own optional parameters, such as a maximum number of iterations.

As an illustration, MOVS uses two single-objective solvers in the two steps of an iteration. The first step is the vertex selection step, which requires a quadratic solver. The second step is the scalarization step, which requires another solver. To use GUROBI for the vertex selection step, with verbosity activated and MOSEK for the second step with a tolerance equal to \(10^{-4}\), within \(140\) iterations, one writes:
\begin{minted}[frame=lines]{python}
pb.solve(solver="MOVS",
         max_iter=140,
         scalarization_solver_options={"solver": cp.MOSEK, "eps": 1e-4},
         vertex_selection_solver_options={"solver": cp.GUROBI, "verbose": True})
\end{minted}

\subsection{Some implementation details}

This section is dedicated to implementation details that are rarely mentioned in generic descriptions of scalarization-based algorithms (but are still present in their implementations). Note that this is not a bad thing. Adding such details would complicate the procedure's description and is irrelevant to its mathematical analysis. However, these details may be pertinent for users who would want to implement their own method.  

\subsubsection{Handling failures of single-objective subproblems}

The three algorithms implemented in MOCVXPY all have an initialization phase that uses the weighted-sum scalarization. This formulation retains the same decision space description as the original \eqref{ref:MOP}/\eqref{ref:VOP} problem. Consequently, it is possible to detect whether the problem is infeasible or unbounded ``from below'' before entering the main step. In this case, MOCVXPY stops and returns a corresponding \texttt{unbounded} or \texttt{infeasible} flag. A problem can also be unbounded ``from above'' when one or more objective components of a solution generated by the algorithm are equal to \(+ \infty\). In MOCVXPY, a solution's objective component is marked as infinite if its value exceeds \(10^{13}\). MOCVXPY stops the algorithm and returns an \texttt{unbounded} flag when this occurs.

A subproblem may fail during the main loop step if the search region is empty, too small, or numerically unstable. The single-objective solver cannot handle this situation or makes insufficient progress toward a weakly efficient solution. As other available implementations/algorithms, e.g., AdeNa,  HyPad or \texttt{MultiObjectiveAlgorithms.jl}, the algorithms implemented in MOCVXPY continue the optimization process until no single-objective subproblems remain to be solved. To avoid revisiting the same subproblems, unsuccessful parameter combinations are saved and systematically compared against before solving a subproblem. The dimensions of these parameters depend on the number of objectives in the original problem multiplied by a small factor (less than \(5\)). For many applications, the number of objectives is small, and the number of subproblems to solve does not exceed the thousands. Storing and comparing them has negligible memory and computational costs. The idea of storing previously explored parameter combinations for scalarization subproblems in a cache can be found for example in~\cite{Ararat2023, Ararat2022, Ararat2024}.

\subsubsection{Exploiting previous optimizations of single-objective subproblems}

For an efficient implementation, it is important to avoid rebuilding the subproblem structure from scratch each time. Only the parameter values of the scalarization formulations change during the main step. MOCVXPY uses the disciplined parameter programming features of CVXPY\footnote{See \url{https://www.cvxpy.org/tutorial/dpp/index.html}.} to reuse the same subproblem formulation throughout the iterations and potentially apply a warm start, which speeds up the resolution of successive subproblems.

\subsubsection{Computing the outer approximation with vertex enumeration}

Both MOVS and MONMO use an outer polyhedral approximation to approximate the upper image, which is refined by adding supporting halfspaces during iterations. The parameters involved in the scalarization formulations of these algorithms are computed using the vertices of the outer approximation. Accurately computing the V-representation of this polyhedron from its H-representation is critical for the performance of MOVS and MONMO. MOCVXPY uses \texttt{cddlib} for this task. Although other polyhedra manipulation libraries exist, \texttt{cddlib} is used because it offers exact (using the GMP library) and floating-point implementations of vertex enumeration, is maintained and has Python bindings. MOCVXPY uses the following strategy to compute the V-representation of the outer approximation. First, it tries to remove redundant inequalities that define the outer approximation using the floating-point procedure of \texttt{cddlib}. Then it uses the exact procedure on the remaining inequalities to obtain the set of vertices of the polyhedron. The components of these vertices are then reconverted into floating-point numbers. The first step removes constraints that are too close numerically to each other and would not be detected using the exact procedure. The second step minimizes the risk of crashing \texttt{cddlib}, which could occur if the floating-point procedure were used, at the detriment of slower computation time.

\subsubsection{Scaled stopping test}

Setting a low threshold value often results in numerical instability: the single-objective solver cannot make sufficient progress on the scalarization subproblems because its feasible space is small. This impedes any improvement in the refinement of the outer approximation. For this reason, the thresholds used in Benson-type or enclosure-based approaches range from \(10^{-1}\) to \(10^{-3}\) (e.g.,~\cite{Ararat2022, Ararat2024, Ararat2023, Dorfler2022, Eichfelder2023a, Eichfelder2024}), compared with threshold values commonly used in single-objective method implementations (from \(10^{-6}\) to \(10^{-8}\)).

Ideally, the stopping test should be invariant under scaling of the objectives. However, when the objective values of all the feasible solutions are close to zero initially, this feature is undesirable. Inspired by~\cite{Waltz2006}, we propose the following stopping tests, which are given in Table~\ref{tab:scaled_stopping_criteria}.

\begin{table}[!ht]
    \centering
    \begin{tabular}{p{1.6cm}p{5.6cm}p{7.75cm}}
    \hline
    Algorithm & Initial sets & Stopping test \\
    \hline \hline
    AdEnA & Local lower bound set \(L^0 = \{\bm{l^{0}}\}\) \newline Local upper bound set \(U^0 = \{\bm{u}^0\}\) & \(w\left(\mathcal{A}\left(L^k, U^k\right)\right) \leq \varepsilon^{opt} \max\{1, w\left(\mathcal{A}\left(L^0, U^0\right)\right)\}\)\\
    \hline
    MONMO & Outer approximation \(\mathcal{O}^0\) & \(d_H(\mathcal{O}^k, \mathcal{P}) \leq \varepsilon^{opt} \max\{1, d_H(\mathcal{O}^0, \mathcal{P})\}\)\\
    \hline
    MOVS & Outer approximation \(\mathcal{O}^0\) \newline Inner approximation \(\mathcal{I}^0\) with \newline \(\mathcal{I}^0 = \operatorname{conv}(f(\mathcal{X}^0)) + \mathcal{C}\) & \(d_H(\mathcal{O}^k, \mathcal{I}^k) \leq \varepsilon^{opt} \max\{1, d_H(\mathcal{O}^0, \mathcal{I}^0)\}\)\\
    \hline
    \end{tabular}
    \caption{Scaled stopping criteria for algorithm implementations in MOCVXPY. The stopping tolerance \(\varepsilon^{opt} > 0\) is provided by the user.}
    \label{tab:scaled_stopping_criteria}
\end{table}

The proposed scaled stopping criteria closely follow those in the original algorithm descriptions. The initial polyhedral outer approximation \(\mathcal{O}^0\), and the local lower bound set \(L^0\) are computed by solving weighted sum scalarization problems (\ref{ref:WS}). The weight vectors correspond to the extreme rays of the dual of the order cone (which comes back to minimizing each component of the objective function separately in the multiobjective case). The AdEnA implementation uses a heuristic to provide the initial local upper bound set \(U^0\). Given a finite set of initial solutions \(\mathcal{X}^0\), used to compute the initial local lower bound set \(L^0\), each coordinate of \(\bm{u}^0\) is defined as \(\bm{u}^0_i = \max\left\{f_i(\bm{x}) + \delta^{\text{AdEnA}} : \bm{x} \in \mathcal{X}^0\right\}\), where \(\delta^{\text{AdEnA}} > 0\) is a positive tolerance value near \(0\). When the maximum is not reached by the quantity metric used by the different algorithms, the scaling of the objectives is taken into account. The factor of \(1\) is used as a safeguard, when the quantity metric is near zero.

\subsection{Parallelism}

To take advantage of multicore processors available on modern machines and speed up problem resolution, a practitioner can select the number of threads the single-objective solver may use to solve a subproblem via the \texttt{scalarization\_solver\_options} argument, as illustrated in Section~\ref{sec:scalarization_formulations}. MOCVXPY also allows one to solve several subproblems in parallel for the MOVS and MONMO methods. The following code shows how to use this feature.
\begin{minted}[frame=lines, linenos]{python}
# Description of the problem
# ...
from dask.distributed import Client

# Define a client that connects to a cluster.
# You can control the number of threads, the number of workers, ...
client = Client(n_workers=2, threads_per_worker=4)
pb.solve(client=client,
         solver="MONMO",
         # Add other options if needed
)
\end{minted}
MOCVXPY uses the Dask library~\cite{Rocklin2015} to manage parallel processing options. This library implements various parallel features that can be deployed and monitored on a single machine or cluster. A \texttt{Client} instance, defined on line 7, controls the number of resources. This instance is then passed to the \texttt{.solve()} method on line \(8\).

\begin{figure}[!th]
    \centering
    \subfigure[Sequential MOVS]{
    \includegraphics[scale=0.5]{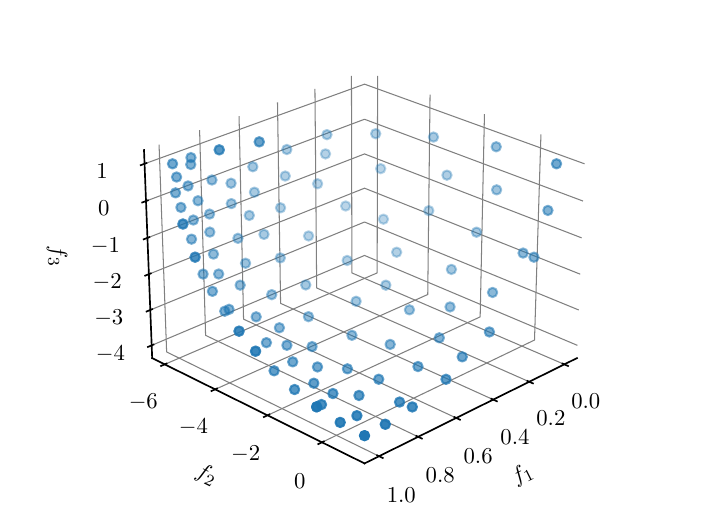}
    }
    \quad
    \subfigure[Parallel MOVS (12 threads)]{
    \includegraphics[scale=0.5]{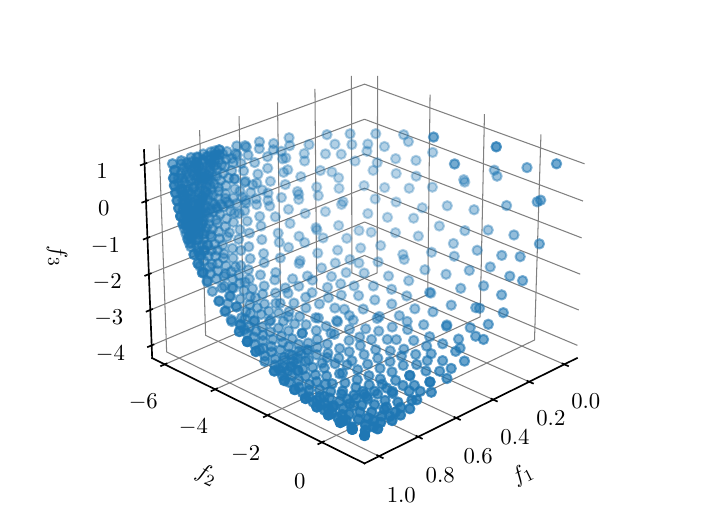}
    }
    \caption{Pareto front approximations using sequential and parallel (with 13 threads) MOVS of the following three-objective problem, taken from~\cite{Dorfler2022}:
    \(\displaystyle\min_{\bm{x} \in \mathbb{R}^3} f(\bm{x}) = \bm{x}\) with respect to \(\leq_{\mathbb{R}^3_+}\), subject to \((x_1 - 1)^2 +  \left(\dfrac{x_2 - 1}{7}\right)^2 + \left(\dfrac{x_3 - 1}{5}\right)^2 \leq 1\).}
    \label{fig:ellipsoid_mo_pb}
\end{figure}

Behind the scene, the \texttt{mocvxpy.problem} instance calls a parallel variant of the given algorithm. At each iteration, this algorithm creates several subproblems instances based on the same scalarization formulations used by the sequential algorithms. The algorithm processes these instances in batches in a parallel loop, without a lock mechanism. This strategy is less thread efficient, since some threads may be idle, but it guarantees reproducible optimization results when run several times with the same parallel options. However, the solution set obtained at the end of the resolution may differ with respect to the corresponding sequential algorithm, since the generation of the scalarization subproblems may change. Figure~\ref{fig:ellipsoid_mo_pb} shows the differences between the MOVS sequential variant and the MOVS parallel variant using 12 threads on the three objective minimization problem. The parallel version obtains more points than the sequential version as it has been designed to exploit more efficiently the parallel ressources it has at its disposal. As a general recommendation, we suggest activating parallelism for large MO/VO problems.

\section{Two real-world applications}\label{section:Real_world_applications}

This section presents two real-world applications that align with the MOCVXPY framework. Additional real-world applications can be found in~\cite{Dorfler2022}.

Their implementations are available and can be found in the \texttt{examples} folder at the root of the project. All experiments described in this section have been conducted on an Apple Mac with an Apple M3 Pro processor with sequential algorithms.

\subsection{Portfolio allocation}

This example, based on the financial Portfolio Selection Problem~\cite{Cornuejols2018OptimizationMethodsFinance, Markowitz1952}, is adapted from a case study from the \texttt{pymoo} website\footnote{\url{https://pymoo.org/case_studies/portfolio_allocation.html}}. The original case study uses a multiobjective evolutionary algorithm to solve this problem.

The goal is to allocate capital to maximize the expected returns of the obtained portfolio and minimize the risk (quantified by an expected variance). Given \(d > 1\) assets, the investor must decide how to split his capital between the different assets, i.e., to design an optimal allocation strategy \(\bm{x} \in \mathbb{R}^d\), where \(\bm{x}\) is a percentage vector. Let \(\bm{r} \in \mathbb{R}^d\) be the expected return of the \(d\) assets: the \(i\)-th coordinate \(r_i\) of \(\bm{r}\) corresponds to the expected return of the \(i\)-th asset for \(i = 1, 2, \ldots, d\). Let \(Q \in \mathbb{R}^{d \times d}\) be the covariance matrix of asset returns. \(Q\) is a positive definite symmetric matrix. The problem is given as
\begin{equation*}
\begin{array}{ll}
    \displaystyle\min_{\bm{x} \in \mathbb{R}^d} & \left(-\bm{r}^\top \bm{x}, \hspace{0.2em} \sqrt{\bm{x}^\top Q \bm{x}}\right)^\top \\
    \text{subject to} & \displaystyle\sum_{i = 1}^d x_i = 1 \\
    & \bm{x} \geq \bm{0}.
\end{array}
\end{equation*}
The first objective is the expected return of the portfolio to maximize (which is equivalent to minimizing its opposite value) and the second objective represents the expected risk of the portfolio to minimize. The constraints force \(\bm{x}\) to be a percentage. The problem can be alternatively reformulated as a convex square-root free model:
\begin{equation*}
\begin{array}{ll}
    \displaystyle\min_{\bm{x} \in \mathbb{R}^d} & \left(-\bm{r}^\top \bm{x}, \hspace{0.2em}  \bm{x}^\top Q \bm{x}\right)^\top \\
    \text{subject to} & \displaystyle\sum_{i = 1}^d x_i = 1 \\
    & \bm{x} \geq \bm{0}.
\end{array}
\end{equation*}

As an illustration, we use the same dataset as in the \texttt{pymoo} case study to solve this problem. The dataset contains historical data for \(d = 20\) technological companies assets from December 1989 to April 2018. From this dataset, we compute the expected return vector \(\bm{r}\) and the correlation matrix \(Q\). After modeling it with MOCVXPY (default values), we obtain the optimal set of solutions in less than \(2\) seconds.

\begin{figure}[!th]
    \centering
    \includegraphics[width=0.5\linewidth]{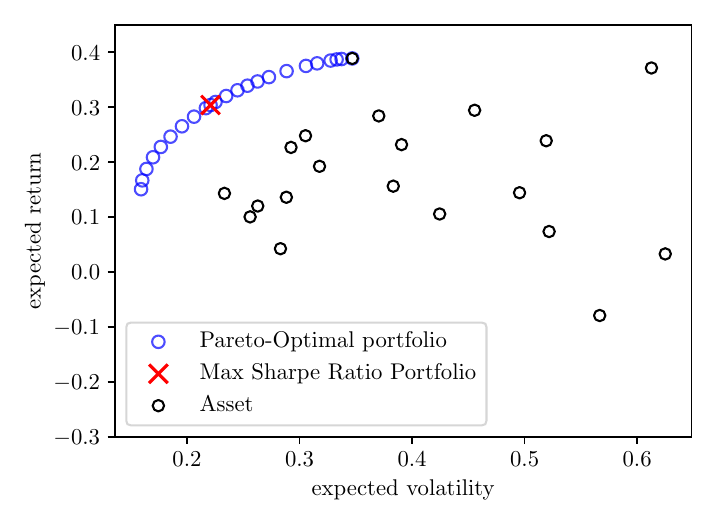}
    \caption{Pareto optimal portfolios and max Sharpe Ratio index portfolio.}
    \label{fig:portfolio_allocation}
\end{figure}

Figure~\ref{fig:portfolio_allocation} shows the optimal set of solutions in the objective space (in blue) as well as the expected risk and expected return of each individual asset (in black). The investor can then see the different optimal trade-offs between the objectives when selecting an optimal strategy. The strategy that maximizes the Sharpe Ratio index~\cite{Cornuejols2018OptimizationMethodsFinance}, commonly used in Portfolio allocation, is indicated in red.

\subsection{Multiobjective optimal powerflow of an electric network via semidefinite relaxation}

This example is inspired by the two following works~\cite{Davoodi2021, Lupien2025}.

The single-objective alternative current optimal powerflow (AC-OPF) problem consists in minimizing the generation costs of an electric network while satisfying load requirements and transmission constraints, which depend on the network's characteristics~\cite{Taylor2015ConvexOptimization}. In the multiobjective case, one additionally aims at minimizing the power losses in the network (e.g.,~\cite{Davoodi2021}). The description of the mathematical model and notations are adapted from~\cite{Lupien2025}.

The electric power system is modeled as a graph (\(\mathcal{B}, \mathcal{L}\)) where \(\mathcal{B} \subset \mathbb{N}\) is the set of vertices, i.e., buses, and \(\mathcal{L} \subset \mathbb{N} \times \mathbb{N}\) is the set of edges, i.e., powerlines. Let \(\mathcal{G} \subseteq \mathcal{B}\) be the set of generators. Each line \((i, j) \in \mathcal{L}\) has its own admittance, denoted by \(y_{i,j} \in \mathbb{C}\). Each bus \(i \in \mathcal{B}\) is characterized by its active power injection \(p_i \in \mathbb{R}\), its reactive power absorption \(q_i \in \mathbb{R}\), and its complex voltage phasor \(v_i \in \mathbb{C}\). Let \(p_{i,j} \in \mathbb{R}\) and \(q_{i, j} \in \mathbb{R}\) be the respective active and reactive power flow in line \((i, j) \in \mathcal{L}\). Finally, for a complex number \(a \in \mathbb{C}\), its conjugate is denoted as \(a^* \in \mathbb{C}\). The multiobjective AC-OPF considered in this work is given as:
\begin{subequations}
\begin{alignat}{3}
    & \min \quad && \left(\displaystyle\sum_{i \in \mathcal{G}} c_i(p_i), \displaystyle\sum_{i \in \mathcal{G}} p_i\right)^\top && \label{acopf:objectives} \\
    & \text{subject to} \quad && p_{i,j} + \jmath q_{i,j} =  v_i (v_i^* - v_j^*) \ y_{i,j}^* && \quad \text{for } (i,j) \in \mathcal{L}, \label{acopf:eq_powerflow}\\
    &&& v_i^{\min} \leq |v_i| \leq v_i^{\max} && \quad \text{for } i \in \mathcal{B}, \label{acopf:voltage_limits} \\
    &&& p_i = \displaystyle\sum_{(i,j) \in \mathcal{L}} p_{i,j}, \ q_i =  \displaystyle\sum_{(i,j) \in \mathcal{L}} q_{i,j} && \quad \text{for } i \in \mathcal{B}, \label{acopt:nodal_power_balance}\\
    &&& p_i^{\min} \leq p_i \leq p_i^{\max}, \ q_i^{\min} \leq q_i \leq q_i^{\max} && \quad \text{for } i \in \mathcal{B}, \label{acopf:power_limits}\\
    &&& p_{i, j}^2 + q_{i, j}^2 \leq (s^{\max}_{i, j})^2 && \quad \text{for } (i, j) \in \mathcal{L}. \label{acopf:powerline_limits}
\end{alignat}
\end{subequations}
Equation~\eqref{acopf:objectives} describes the two objectives. The first objective represents the generation costs, defined as a convex cost function, where each \(c_i : x \in \mathbb{R} \rightarrow \mathbb{R}\) is a convex quadratic function. The second objective function denotes the active power losses of the electric network, that must also be minimized. Equation~\eqref{acopf:eq_powerflow} is the powerflow constraints. Note that it is nonlinear and non-convex. Equations~\eqref{acopf:voltage_limits},~\eqref{acopf:power_limits} and~\eqref{acopf:powerline_limits} respectively impose voltage magnitude limits, active and reactive power limits on buses and apparent power limit \(s_{i,j}^{\max}\) on line \((i,j) \in \mathcal{L}\). Equation~\eqref{acopt:nodal_power_balance} describes the active and reactive nodal power balance.

One can rewrite Equations~\eqref{acopf:eq_powerflow} and~\eqref{acopf:voltage_limits} by lifting them into higher dimensions~\cite{Taylor2015ConvexOptimization}:
\begin{equation}
p_{i,j} + \jmath q_{i,j} = (W_{i,i} - W_{i,j}) y^*_{i, j} \text{ for } (i,j) \in \mathcal{L}
\end{equation}
and
\begin{equation}
\left(v_i^{\min}\right)^2 \leq W_{i,i} \leq \left(v_i^{\max}\right)^2 \text{ for } i \in \mathcal{B}
\end{equation}
where \(W \in \mathbb{C}^{\mathcal{|B| \times |\mathcal{B}|}}\) is a Hermitian matrix. In addition, \(W\) must satisfy:
\begin{subequations}
    \begin{flalign}
        W & \succeq 0, \\
        \text{rank}(W) & = 1. \label{acopf_relax:rank1_cstr}
    \end{flalign}
\end{subequations}
When omitting Equation~\eqref{acopf_relax:rank1_cstr}, one obtains the semi-definite relaxation of (AC-OPF), denoted as (AC-OPF-SDR)~\cite{Taylor2015ConvexOptimization}, which is convex. There exists several conditions related to the topology of the network under which one can prove that this relaxation is exact (e.g., see~\cite{Lupien2025, Taylor2015ConvexOptimization}).

This work uses the electric network described in the case 9\footnote{\url{https://matpower.org/docs/ref/matpower5.0/case9.html}} of MATPOWER~\cite{Zimmerman2011}. Using MOCVXPY, we model the AC-OPF-SDR relaxation problem, resulting in a problem with \(135\) decision variables. The problem is solved with the AdEnA algorithm. MOSEK is selected as the main single-objective solver for the subproblems. Indeed, open-source conic solvers such as SCS or CLARABEL provided with CVXPY/MOCVXPY fail to solve the subproblems, hence this choice.

\begin{figure}
    \centering
    \includegraphics[scale=0.8]{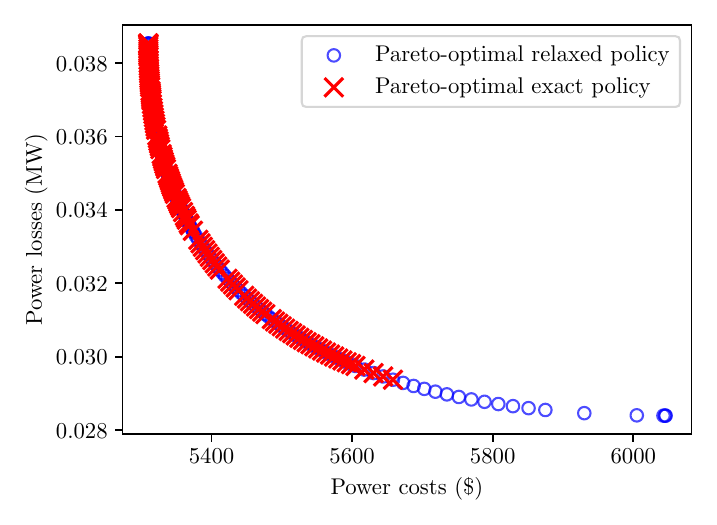}
    \caption{Pareto optimal relaxed and exact solutions for the AC-OPF problem on the case \(9\) of MATPOWER.}
    \label{fig:acopf_ieee9}
\end{figure}

AdEnA takes less than 2 seconds to solve the AC-OPF-SDR, i.e., to reach the stopping tolerance. Following~\cite{Lupien2025}, a postprocessing step is performed to check if the obtained solutions are exact. Figure~\ref{fig:acopf_ieee9} displays the optimal solutions obtained at the end of the resolution for the AC-OPF-SDR relaxation. Among them, the red-crossed ones are also exact for AC-OPF. As can be seen, most of them are exact.

\begin{remark}
It is important to emphasize that solving the SDR relaxation to obtain exact solutions of AC-OPF may not work for all electric networks. Among all the instances described in MATPOWER, the case 9 is a ``small electric'' network, which has been chosen to illustrate the potential of MOCVXPY.
\end{remark}

\section{Discussion}\label{section:Discussion}

As emphasized in a recent survey on exact multiobjective optimization methods, ``it is important that the [new] developed algorithms are provided to the public such that multiobjective optimization techniques are used more widely for solving application problems in practice''~\cite{Eichfelder2021b}. In the lineage of these recommendations, this paper introduces MOCVXPY, an open-source Python library, to solve convex multiobjective and vector optimization problems.

Built on top of CVXPY, it reuses most of its API to allow the user to describe his/her model in an intuitive way that follows its mathematical formulation and solve it. Currently, it provides three state-of-the-art algorithms possessing convergence guarantees for convex multiobjective/vector optimization. In addition, the practitioner has access to basic parallel functionalities if his/her problem is large. Two applications that fit into the framework of MOCVXPY are described and solved.

For the moment, MOCVXPY does not offer methods to solve mixed-integer convex problems. Future work could involve the integration of such methods. All methods implemented into this library assume that the problems are bounded. For these last three years, researchers have started to tackle unbounded convex problems: methods such as~\cite{Wagner2023} could be added to MOCVXPY. Finally, the solvers currently provided remain generic. For example, it is possible to solve MOLP or VOLPs using MOCVXPY, but these algorithms will not always capture the exact ``shape'' of the upper image (which is a polyhedron). Solvers such as \texttt{Bensolve}, \texttt{Inner} or \texttt{polySCIP} would be more relevant. A deeper integration of these solvers for the linear case is planned.

\paragraph{Acknowledgments} The authors would like to thank Professor Gabriele Eichfelder and Dr. Leo Warnow for their valuable explanations on AdEnA. 

\paragraph{Fundings} This work was supported by the Friedrich Schiller University Jena under the Jena Excellence Fellowship Program.

\section*{Declaration}

\paragraph{Conflict of interest} The authors report no conflict of interest.

\bibliographystyle{plain}
\bibliography{bibliography}

\begin{thebibliography}{10}

\bibitem{Alexandropoulos2019}
S.-A.~N. Alexandropoulos, C.~K. Aridas, S.~B. Kotsiantis, and M.~N. Vrahatis.
\newblock {\em Multi-Objective Evolutionary Optimization Algorithms for Machine
  Learning: A Recent Survey}, pages 35--55.
\newblock Springer International Publishing, 2019.

\bibitem{Ararat2023a}
{\c C}.~Ararat and N.~Meimanjan.
\newblock Computation of {{Systemic Risk Measures}}: {{A Mixed-Integer
  Programming Approach}}.
\newblock {\em Operations Research}, 71(6):2130--2145, 2023.

\bibitem{Ararat2023}
{\c C}.~Ararat, S.~Tekg{\"u}l, and F.~Ulus.
\newblock Geometric {{Duality Results}} and {{Approximation Algorithms}} for
  {{Convex Vector Optimization Problems}}.
\newblock {\em SIAM Journal on Optimization}, 33(1):116--146, 2023.

\bibitem{Ararat2022}
{\c C}.~Ararat, F.~Ulus, and M.~Umer.
\newblock A {{Norm Minimization-Based Convex Vector Optimization Algorithm}}.
\newblock {\em Journal of Optimization Theory and Applications},
  194(2):681--712, 2022.

\bibitem{Ararat2024}
{\c C}.~Ararat, F.~Ulus, and M.~Umer.
\newblock Convergence analysis of a norm minimization-based convex vector
  optimization algorithm.
\newblock {\em SIAM Journal on Optimization}, 34(3):2700--2728, 2024.

\bibitem{Benson1998a}
H.~P. Benson.
\newblock An {{Outer Approximation Algorithm}} for {{Generating All Efficient
  Extreme Points}} in the {{Outcome Set}} of a {{Multiple Objective Linear
  Programming Problem}}.
\newblock {\em Journal of Global Optimization}, 13(1):1--24, 1998.

\bibitem{Biscani2020}
F.~Biscani and D.~Izzo.
\newblock A parallel global multiobjective framework for optimization: Pagmo.
\newblock {\em Journal of Open Source Software}, 5(53):2338, 2020.

\bibitem{Blank2020}
J.~Blank and K.~Deb.
\newblock Pymoo: {{Multi-objective}} optimization in python.
\newblock {\em IEEE access : practical innovations, open solutions},
  8:89497--89509, 2020.

\bibitem{Borndoerfer2016}
R.~Bornd{\"o}rfer, S.~Schenker, M.~Skutella, and T.~Strunk.
\newblock Polyscip.
\newblock In G.-M. Greuel, T.~Koch, P.~Paule, and A.~Sommese, editors, {\em
  Mathematical Software - ICMS 2016, 5th International Conference, Berlin,
  Germany, July 11-14, 2016, Proceedings}, volume 9725, pages 259 -- 264, 2016.

\bibitem{Branke2008MultiobjectiveOptimizationInteractiveEvolutionaryApproaches}
J.~Branke, K.~Deb, K.~Miettinen, and R.~Slowi{\'n}ski.
\newblock {\em Multiobjective Optimization: {{Interactive}} and Evolutionary
  Approaches}, volume 5252.
\newblock Springer Science \& Business Media, 2008.

\bibitem{Burachik2017}
R.~S. Burachik, C.~Y. Kaya, and M.~M. Rizvi.
\newblock A {{New Scalarization Technique}} and {{New Algorithms}} to
  {{Generate Pareto Fronts}}.
\newblock {\em SIAM Journal on Optimization}, 27(2):1010--1034, 2017.

\bibitem{Burachik2022}
R.~S. Burachik, C.~Y. Kaya, and M.~M. Rizvi.
\newblock Algorithms for generating {{Pareto}} fronts of multi-objective
  integer and mixed-integer programming problems.
\newblock {\em Engineering Optimization}, 54(8):1413--1425, 2022.

\bibitem{Bussieck2004}
M.~R. Bussieck and A.~Meeraus.
\newblock General algebraic modeling system (gams).
\newblock In {\em Modeling languages in mathematical optimization}, pages
  137--157. Springer, 2004.

\bibitem{Bynum2021Pyomo}
M.~L. Bynum, G.~A. Hackebeil, W.~E. Hart, C.~D. Laird, B.~L. Nicholson, J.~D.
  Siirola, J.-P. Watson, D.~L. Woodruff, et~al.
\newblock {\em Pyomo-optimization modeling in python}, volume~67.
\newblock Springer, 2021.

\bibitem{Ciripoi2019}
D.~Ciripoi, A.~L\"ohne, and B.~Wei\ss~ing.
\newblock Calculus of convex polyhedra and polyhedral convex functions by
  utilizing a multiple objective linear programming solver.
\newblock {\em Optimization}, 68(10):2039--2054, 2019.

\bibitem{Collette2011MultiobjectiveOptimizationPrinciplesCaseStudies}
Y.~Collette and P.~Siarry.
\newblock {\em Multiobjective Optimization: Principles and Case Studies}.
\newblock Springer Science \& Business Media, 2011.

\bibitem{Cornuejols2018OptimizationMethodsFinance}
G.~Cornu{\'e}jols, J.~Pe{\~n}a, and R.~T{\"u}t{\"u}nc{\"u}.
\newblock {\em Optimization Methods in Finance}.
\newblock Cambridge University Press, Cambridge, 2 edition, 2018.

\bibitem{Csirmaz2021}
L.~Csirmaz.
\newblock Inner approximation algorithm for solving linear multiobjective
  optimization problems.
\newblock {\em Optimization}, 70(7):1487--1511, 2021.

\bibitem{Cui2017}
Y.~Cui, Z.~Geng, Q.~Zhu, and Y.~Han.
\newblock Review: {{Multi-objective}} optimization methods and application in
  energy saving.
\newblock {\em Energy}, 125:681--704, 2017.

\bibitem{Dachert2017}
K.~D{\"a}chert, K.~Klamroth, R.~Lacour, and D.~Vanderpooten.
\newblock Efficient computation of the search region in multi-objective
  optimization.
\newblock {\em European Journal of Operational Research}, 260(3):841--855,
  2017.

\bibitem{Das1998}
I.~Das and J.~E. Dennis.
\newblock Normal-boundary intersection: A new method for generating the pareto
  surface in nonlinear multicriteria optimization problems.
\newblock {\em SIAM Journal on Optimization}, 8(3):631--657, 1998.

\bibitem{Davoodi2021}
E.~Davoodi, E.~Babaei, B.~Mohammadi-Ivatloo, M.~Shafie-Khah, and J.~P.~S.
  Catalão.
\newblock Multiobjective optimal power flow using a semidefinite
  programming-based model.
\newblock {\em IEEE Systems Journal}, 15(1):158--169, 2021.

\bibitem{DeSantis2021}
M.~De~Santis and G.~Eichfelder.
\newblock A decision space algorithm for multiobjective convex quadratic
  integer optimization.
\newblock {\em Computers \& Operations Research}, 134:105396, 2021.

\bibitem{DeSantis2020}
M.~De~Santis, G.~Eichfelder, J.~Niebling, and S.~Rockt{\"a}schel.
\newblock Solving multiobjective mixed integer convex optimization problems.
\newblock {\em SIAM Journal on Optimization}, 30(4):3122--3145, 2020.

\bibitem{DeSantis2024}
M.~De~Santis, G.~Eichfelder, D.~Patria, and L.~Warnow.
\newblock Using dual relaxations in multiobjective mixed-integer convex
  quadratic programming.
\newblock {\em Journal of Global Optimization}, 2024.

\bibitem{Diamond2016}
S.~Diamond and S.~Boyd.
\newblock Cvxpy: A python-embedded modeling language for convex optimization.
\newblock {\em Journal of Machine Learning Research}, 17(83):1--5, 2016.

\bibitem{Dorfler2022a}
D.~D{\"o}rfler.
\newblock On the {{Approximation}} of {{Unbounded Convex Sets}} by
  {{Polyhedra}}.
\newblock {\em Journal of Optimization Theory and Applications},
  194(1):265--287, 2022.

\bibitem{Dorfler2022}
D.~D{\"o}rfler, A.~L{\"o}hne, C.~Schneider, and B.~Wei{\ss}ing.
\newblock A {{Benson-type}} algorithm for bounded convex vector optimization
  problems with vertex selection.
\newblock {\em Optimization Methods and Software}, 37(3):1006--1026, 2022.

\bibitem{Dowson2025}
O.~Dowson, X.~Gandibleux, and G.~Kof.
\newblock Multiobjectivealgorithms.jl: a julia package for solving
  multi-objective optimization problems, 2025.

\bibitem{Dunning2017}
I.~Dunning, J.~Huchette, and M.~Lubin.
\newblock Jump: A modeling language for mathematical optimization.
\newblock {\em SIAM review}, 59(2):295--320, 2017.

\bibitem{Durillo2011}
J.~J. Durillo and A.~J. Nebro.
\newblock {{jMetal}}: {{A Java}} framework for multi-objective optimization.
\newblock {\em Advances in Engineering Software}, 42(10):760--771, 2011.

\bibitem{Ehrgott2005MulticriteriaOptimization}
M.~Ehrgott.
\newblock {\em Multicriteria Optimization}, volume 491.
\newblock Springer-Verlag, 2005.

\bibitem{Ehrgott2012}
M.~Ehrgott, A.~L{\"o}hne, and L.~Shao.
\newblock A dual variant of {{Benson}}'s ``outer approximation algorithm'' for
  multiple objective linear programming.
\newblock {\em Journal of Global Optimization}, 52(4):757--778, April 2012.

\bibitem{Ehrgott2011}
M.~Ehrgott, L.~Shao, and A.~Sch{\"o}bel.
\newblock An approximation algorithm for convex multi-objective programming
  problems.
\newblock {\em Journal of Global Optimization}, 50(3):397--416, 2011.

\bibitem{Eichfelder2009}
G.~Eichfelder.
\newblock An adaptive scalarization method in multiobjective optimization.
\newblock {\em SIAM Journal on Optimization}, 19(4):1694--1718, 2009.

\bibitem{Eichfelder2014}
G.~Eichfelder.
\newblock Vector optimization in medical engineering.
\newblock In P.~M. Pardalos and T.~M. Rassias, editors, {\em Mathematics
  Without Boundaries: Surveys in Interdisciplinary Research}, pages 181--215.
  Springer, 2014.

\bibitem{Eichfelder2021b}
G.~Eichfelder.
\newblock Twenty years of continuous multiobjective optimization in the
  twenty-first century.
\newblock {\em EURO Journal on Computational Optimization}, 9:100014, 2021.

\bibitem{Eichfelder2021a}
G.~Eichfelder, P.~Kirst, L.~Meng, and O.~Stein.
\newblock A general branch-and-bound framework for continuous global
  multiobjective optimization.
\newblock {\em Journal of Global Optimization}, 80:195--227, 2021.

\bibitem{Eichfelder2024a}
G.~Eichfelder, O.~Stein, and L.~Warnow.
\newblock A solver for multiobjective mixed-integer convex and nonconvex
  optimization.
\newblock {\em Journal of Optimization Theory and Applications},
  203:1736--1766, 2024.

\bibitem{Eichfelder2025}
G.~Eichfelder and F.~Ulus.
\newblock Local upper bounds based on polyhedral ordering cones.
\newblock Technical report, 2025.

\bibitem{Eichfelder2021}
G.~Eichfelder and L.~Warnow.
\newblock An approximation algorithm for multi-objective optimization problems
  using a box-coverage.
\newblock {\em Journal of Global Optimization}, 83(2):329--357, 2021.

\bibitem{Eichfelder2023a}
G.~Eichfelder and L.~Warnow.
\newblock Advancements in the computation of enclosures for multi-objective
  optimization problems.
\newblock {\em European Journal of Operational Research}, 310(1):315--327,
  2023.

\bibitem{Eichfelder2024}
G.~Eichfelder and L.~Warnow.
\newblock A hybrid patch decomposition approach to compute an enclosure for
  multi-objective mixed-integer convex optimization problems.
\newblock {\em Mathematical Methods of Operations Research}, 100(1):291--320,
  2024.

\bibitem{Feinstein2017}
Z.~Feinstein and B.~Rudloff.
\newblock A recursive algorithm for multivariate risk measures and a set-valued
  {{Bellman}}'s principle.
\newblock {\em Journal of Global Optimization}, 68(1):47--69, 2017.

\bibitem{Feinstein2024a}
Z.~Feinstein and B.~Rudloff.
\newblock Technical {{Note}}---{{Characterizing}} and {{Computing}} the {{Set}}
  of {{Nash Equilibria}} via {{Vector Optimization}}.
\newblock {\em Operations Research}, 72(5):2082--2096, 2024.

\bibitem{Fortin2012}
F.-A. Fortin, F.-M. {De Rainville}, M.-A. Gardner, M.~Parizeau, and C.~Gagn\'e.
\newblock {DEAP}: Evolutionary algorithms made easy.
\newblock {\em Journal of Machine Learning Research}, 13:2171--2175, 2012.

\bibitem{Fourer1990}
R.~Fourer, D.~M. Gay, and B.~W. Kernighan.
\newblock Ampl: A mathematical programming language.
\newblock {\em Management Science}, 36(5):519--554, 1990.

\bibitem{Fukuda2014}
E.~H. Fukuda and L.~M.~G. Drummond.
\newblock A {{Survey}} on {{Multiobjective Descent Methods}}.
\newblock {\em Pesquisa Operacional}, 34:585--620, 2014.

\bibitem{Garett2012}
D.~Garrett.
\newblock inspyred (version 1.0.1 [software]. inspired intelligence, 2012.

\bibitem{Grant2008}
M.~Grant and S.~Boyd.
\newblock Graph implementations for nonsmooth convex programs.
\newblock In V.~Blondel, S.~Boyd, and H.~Kimura, editors, {\em Recent Advances
  in Learning and Control}, Lecture Notes in Control and Information Sciences,
  pages 95--110. Springer-Verlag Limited, 2008.
\newblock \url{http://stanford.edu/~boyd/graph_dcp.html}.

\bibitem{Grant2014}
M.~Grant and S.~Boyd.
\newblock {CVX}: Matlab software for disciplined convex programming, version
  2.1.
\newblock \url{https://cvxr.com/cvx}, 2014.

\bibitem{Greer1984}
R.~Greer.
\newblock Chapter 2: A tutorial on polyhedral convex cones.
\newblock In R.~Greer, editor, {\em Trees and Hills: Methodology for Maximizing
  Functions of Systems of Linear Relations}, volume~96 of {\em North-Holland
  Mathematics Studies}, pages 15--81. North-Holland, 1984.

\bibitem{Hadka2024}
D.~Hadka.
\newblock Platypus: A framework for evolutionary computing in python (version
  1.4.1) [computer software], 2024.

\bibitem{Hamel2013}
A.~H. Hamel, B.~Rudloff, and M.~Yankova.
\newblock Set-valued average value at risk and its computation.
\newblock {\em Mathematics and Financial Economics}, 7(2):229--246, 2013.

\bibitem{Handl2007}
J.~Handl, D.~B. Kell, and J.~Knowles.
\newblock Multiobjective optimization in bioinformatics and computational
  biology.
\newblock {\em IEEE/ACM Transactions on Computational Biology and
  Bioinformatics}, 4(2):279--292, 2007.

\bibitem{Helfrich2024a}
S.~Helfrich, S.~Ruzika, and C.~Thielen.
\newblock Efficiently constructing convex approximation sets in multiobjective
  optimization problems.
\newblock {\em INFORMS Journal on Computing}, 2024.

\bibitem{Heyde2011}
F.~Heyde and A.~L{\"o}hne.
\newblock Solution concepts in vector optimization: A fresh look at an old
  story.
\newblock {\em Optimization}, 60(12):1421--1440, 2011.

\bibitem{Jahn2009VectorOptimization}
J.~Jahn.
\newblock {\em Vector Optimization}.
\newblock Springer, 2009.

\bibitem{Kaibel2011}
V.~Kaibel.
\newblock {\em Basic Polyhedral Theory}, pages 396--409.
\newblock John Wiley \& Sons, Ltd, 2011.

\bibitem{Keskin2023}
{\.I}.~N. Keskin and F.~Ulus.
\newblock Outer approximation algorithms for convex vector optimization
  problems.
\newblock {\em Optimization Methods and Software}, 38(4):723--755, 2023.

\bibitem{Klamroth2015}
K.~Klamroth, R.~Lacour, and D.~Vanderpooten.
\newblock On the representation of the search region in multi-objective
  optimization.
\newblock {\em European Journal of Operational Research}, 245(3):767--778,
  2015.

\bibitem{Lammel2025}
I.~Lammel, K.-H. K{\"u}fer, and P.~S{\"u}ss.
\newblock Efficient {{Approximation Quality Computation}} for {{Sandwiching
  Algorithms}} for {{Convex Multicriteria Optimization}}.
\newblock {\em Journal of Optimization Theory and Applications}, 204(3):41,
  2025.

\bibitem{Lofberg2004}
J.~Lofberg.
\newblock Yalmip : a toolbox for modeling and optimization in matlab.
\newblock In {\em 2004 IEEE International Conference on Robotics and Automation
  (IEEE Cat. No.04CH37508)}, pages 284--289, 2004.

\bibitem{Lohne2011Vector}
A.~L{\"o}hne.
\newblock {\em Vector optimization with infimum and supremum}.
\newblock Springer Science \& Business Media, 2011.

\bibitem{Lohne2014a}
A.~L{\"o}hne and B.~Rudloff.
\newblock An algorithm for calculating the set of superhedging portfolios in
  markets with transaction costs.
\newblock {\em International Journal of Theoretical and Applied Finance},
  17(02):1450012, 2014.

\bibitem{Lohne2014}
A.~L{\"o}hne, B.~Rudloff, and F.~Ulus.
\newblock Primal and dual approximation algorithms for convex vector
  optimization problems.
\newblock {\em Journal of Global Optimization}, 60(4):713--736, 2014.

\bibitem{Lohne2017}
A.~L{\"o}hne and B.~Wei{\ss}ing.
\newblock The vector linear program solver {{{\emph{Bensolve}}}} -- notes on
  theoretical background.
\newblock {\em European Journal of Operational Research}, 260(3):807--813,
  2017.

\bibitem{Lupien2025}
J.-L. Lupien and A.~Lesage-Landry.
\newblock Ex post conditions for the exactness of optimal power flow conic
  relaxations.
\newblock {\em Electric Power Systems Research}, 238:111130, 2025.

\bibitem{Markowitz1952}
H.~Markowitz.
\newblock Portfolio selection.
\newblock {\em Journal of Finance}, 7(1):77--91, 1952.

\bibitem{Miettinen1999NonlinearMultiobjectiveOptimization}
K.~Miettinen.
\newblock {\em Nonlinear Multiobjective Optimization}, volume~12.
\newblock Springer Science \& Business Media, 1999.

\bibitem{Niebling2019}
J.~Niebling and G.~Eichfelder.
\newblock A branch--and--bound-based algorithm for nonconvex multiobjective
  optimization.
\newblock {\em SIAM Journal on Optimization}, 29(1):794--821, January 2019.

\bibitem{Raimundo2020}
M.~M. Raimundo, P.~A.~V. Ferreira, and F.~J. Von~Zuben.
\newblock An extension of the non-inferior set estimation algorithm for many
  objectives.
\newblock {\em European Journal of Operational Research}, 284(1):53--66, 2020.

\bibitem{Rocklin2015}
M.~Rocklin et~al.
\newblock Dask: Parallel computation with blocked algorithms and task
  scheduling.
\newblock In {\em Proceedings of the 14th python in science conference
  (SciPy)}, pages 126--132, 2015.

\bibitem{Rudloff2021}
B.~Rudloff and F.~Ulus.
\newblock Certainty equivalent and utility indifference pricing for incomplete
  preferences via convex vector optimization.
\newblock {\em Mathematics and Financial Economics}, 15(2):397--430, March
  2021.

\bibitem{Ruszczynski2003}
A.~Ruszczynski and R.~J. Vanderbei.
\newblock Frontiers of {{Stochastically Nondominated Portfolios}}.
\newblock {\em Econometrica}, 71(4):1287--1297, 2003.

\bibitem{Sagnol2022}
G.~Sagnol and M.~Stahlberg.
\newblock Picos: A python interface to conic optimization solvers.
\newblock {\em Journal of Open Source Software}, 7(70):3915, 2022.

\bibitem{Schutze2016a}
O.~Sch{\"u}tze, A.~Mart{\'i}n, A.~Lara, S.~Alvarado, E.~Salinas, and
  C.~A.~Coello Coello.
\newblock The directed search method for multi-objective memetic algorithms.
\newblock {\em Computational Optimization and Applications}, 63(2):305--332,
  2016.

\bibitem{Shami2022}
T.~M. Shami, A.~A. {El-Saleh}, M.~Alswaitti, Q.~{Al-Tashi}, M.~A. Summakieh,
  and S.~Mirjalili.
\newblock Particle {{Swarm Optimization}}: {{A Comprehensive Survey}}.
\newblock {\em IEEE Access}, 10:10031--10061, 2022.

\bibitem{Sharma2022}
S.~Sharma and V.~Kumar.
\newblock A {{Comprehensive Review}} on {{Multi-objective Optimization
  Techniques}}: {{Past}}, {{Present}} and {{Future}}.
\newblock {\em Archives of Computational Methods in Engineering},
  29(7):5605--5633, 2022.

\bibitem{Sharma2013}
S.~Sharma and G.~P. Rangaiah.
\newblock Multi-objective optimization applications in chemical engineering.
\newblock In {\em Multi-objective Optimization in Chemical Engineering},
  chapter~3, pages 35--102. John Wiley \& Sons, Ltd, 2013.

\bibitem{Taylor2015ConvexOptimization}
J.~A. Taylor.
\newblock {\em Convex Optimization of Power Systems}.
\newblock Cambridge University Press, 2015.

\bibitem{Wagner2023}
A.~Wagner, F.~Ulus, B.~Rudloff, G.~Kov{\'a}{\v c}ov{\'a}, and N.~Hey.
\newblock Algorithms to {{Solve Unbounded Convex Vector Optimization
  Problems}}.
\newblock {\em SIAM Journal on Optimization}, 33(4):2598--2624, 2023.

\bibitem{Waltz2006}
R.~A. Waltz, J.~L. Morales, J.~Nocedal, and D.~Orban.
\newblock An interior algorithm for nonlinear optimization that combines line
  search and trust region steps.
\newblock {\em Mathematical programming}, 107(3):391--408, 2006.

\bibitem{Wiecek2016}
M.~M. Wiecek, M.~Ehrgott, and A.~Engau.
\newblock Continuous multiobjective programming.
\newblock In {\em Multiple {{Criteria Decision Analysis}}}, pages 739--815.
  Springer New York, 2016.

\bibitem{Zimmerman2011}
R.~D. Zimmerman, C.~E. Murillo-Sánchez, and R.~J. Thomas.
\newblock {MATPOWER}: Steady-state operations, planning, and analysis tools for
  power systems research and education.
\newblock {\em IEEE Transactions on Power Systems}, 26(1):12--19, 2011.

\end{thebibliography}

\end{document}